\documentclass[journal,twoside]{IEEEtran}

\usepackage{cite}
\usepackage[cmex10]{amsmath,mathtools}
\usepackage{amsfonts}
\usepackage{amssymb}
\usepackage{textcomp}
\usepackage[caption=false,font=footnotesize]{subfig}
\usepackage{fixltx2e}
\usepackage{url}
\usepackage{comment}
\usepackage{siunitx}
\usepackage{multirow}
\usepackage{tikz}
\usetikzlibrary{arrows,decorations.markings,positioning}
\usepackage[european]{circuitikz}
\usepackage{pgfplots}
\usepackage{xspace}
\usepackage{accents}
\usepackage{ifthen}
\usepackage{bm}

\MakeRobust{\underaccent}

\DeclareSymbolFont{largesymbols}{OMX}{yhex}{m}{n}
\DeclareMathAccent{\wideparen}{\mathord}{largesymbols}{"F3}

\newtheorem{theorem}{Theorem}
\newtheorem{lemma}{Lemma}
\newtheorem{proposition}{Proposition}
\newtheorem{corollary}{Corollary}
\newtheorem{definition}{Definition}
\newtheorem{requirement}[definition]{Definition}
\newtheorem{assumption}{Assumption}
\newtheorem{definitionx}[definition]{Definition}
\newtheorem{remark}{Remark}

\newcommand{\tabstretch}{1.4}
\newcommand{\costnum}[1]{\$$\num[group-separator={,}]{#1}\,$/\,h}

\DeclareMathOperator{\real}{Re}
\DeclareMathOperator{\imag}{Im}
\DeclareMathOperator{\trace}{tr}
\DeclareMathOperator{\rank}{rank}

\DeclareMathOperator*{\minimize}{minimize}
\DeclareMathOperator*{\subjectto}{subject~to}

\DeclareMathOperator{\cone}{cone}

\DeclareMathOperator{\interior}{int}

\DeclareMathOperator{\psd}{\succeq}
\DeclareMathOperator{\pd}{\succ}
\DeclareMathOperator{\en}{e}

\newcommand{\abs}[1]{{\lvert#1\rvert}}

\newcommand{\vm}[1]{{\bm{#1}}}
\newcommand{\mc}[1]{{\mathcal{#1}}}

\newcommand{\conj}{^\ast}
\DeclareMathOperator{\transpose}{T}
\DeclareMathOperator{\hermitian}{H}
\newcommand{\tran}{^{\transpose}}
\newcommand{\herm}{^{\hermitian}}
\newcommand{\nR}{\mathbb{R}}
\newcommand{\nRnn}{\mathbb{R}_{+}}

\newcommand{\nC}{\mathbb{C}}
\newcommand{\nS}{\mathbb{S}}
\newcommand{\aOut}[1]{\expandafter\hat#1}
\newcommand{\aIn}[1]{\expandafter\check#1}
\newcommand{\aBar}[1]{\expandafter\bar#1}
\newcommand{\aUBar}[1]{\expandafter\underaccent{\bar}#1}
\newcommand{\aShunt}[1]{\expandafter\tilde#1}
\newcommand{\aArrow}[1]{\expandafter\vec#1}
\newcommand{\aDot}[1]{\expandafter\dot#1}
\newcommand{\aLoss}[1]{\expandafter\tilde#1}
\newcommand{\aCost}[1]{\expandafter\bar#1}
\newcommand{\aDc}[1]{\expandafter\mathring#1}

\newcommand{\aPgen}[1]{#1^{\text{(G)}}}
\newcommand{\aPload}[1]{#1^{\text{(L)}}}
\newcommand{\aPgenmin}[1]{\expandafter\underaccent{\bar}#1^{\text{(G)}}}
\newcommand{\aPgenmax}[1]{\expandafter\bar#1^{\text{(G)}}}
\newcommand{\aLB}[1]{\expandafter\underaccent{\bar}#1}
\newcommand{\aUB}[1]{\expandafter\bar#1}

\newcommand{\aOutUB}[1]{\bar{\hat{#1}}}

\newcommand{\aInUB}[1]{\bar{\check{#1}}}

\newcommand{\aDcLB}[1]{\underaccent{\bar}{\mathring{#1}}}
\newcommand{\aDcUB}[1]{\bar{\mathring{#1}}}
\newcommand{\aOpt}[1]{#1^{\star}}

\newcommand{\aArc}[1]{\expandafter\wideparen#1}

\newcommand{\sV}{{\mc{V}}}
\newcommand{\sE}{{\mc{E}}}
\newcommand{\sD}{{\mc{D}}}
\newcommand{\sB}{{\mc{B}}}
\newcommand{\sBE}{\sB_\sE}
\newcommand{\sBD}{\sB_\sD}

\newcommand{\sS}{{\mc{S}}}
\newcommand{\sC}{{\mc{C}}}
\newcommand{\sH}{{\mc{H}}}
\newcommand{\sEU}{\aBar{\sE}}
\newcommand{\sX}{{\mc{X}}}

\newcommand{\gG}{{\mc{G}}}
\newcommand{\faE}{{\aOut{\epsilon}}}
\newcommand{\fbE}{{\aIn{\epsilon}}}
\newcommand{\faD}{{\aOut{\delta}}}
\newcommand{\fbD}{{\aIn{\delta}}}
\newcommand{\fPsi}{\varPsi}
\newcommand{\fpsi}{\psi}
\newcommand{\fPsiCone}{\hat{\fPsi}}
\newcommand{\fpsiCone}{\hat{\fpsi}}
\newcommand{\fMtxGraph}[1]{{\mc{T}}(#1)}
\newcommand{\cNV}{{N}_\sV}
\newcommand{\cNE}{{N}_\sE}
\newcommand{\cND}{{N}_\sD}
\newcommand{\cPF}{{\text{PF}}}
\newcommand{\mY}{{\vm{Y}}}

\newcommand{\mM}{{\vm{M}}}
\newcommand{\mV}{{\vm{V}}}
\newcommand{\mI}{{\vm{I}}}
\newcommand{\mIout}{{\hat{\skew{-3}\bm{I}}\hskip-2.25pt}}
\newcommand{\mIin}{{\check{\skew{-3}\bm{I}}\hskip-2.25pt}}
\newcommand{\mS}{{\vm{S}}}
\newcommand{\mP}{{\vm{P}}}
\newcommand{\mQ}{{\vm{Q}}}
\newcommand{\mA}{{\vm{A}}}
\newcommand{\mAub}{{\bar{\skew{-6}\bm{A}}\hskip-3.5pt}}
\newcommand{\mB}{{\vm{B}}}
\newcommand{\mC}{{\vm{C}}}

\newcommand{\mZero}{{\vm{0}}}
\newcommand{\vv}{{\vm{v}}}
\newcommand{\vi}{{\vm{i}}}
\newcommand{\ve}{{\vm{e}}}

\newcommand{\vp}{{\vm{p}}}
\newcommand{\vh}{{\vm{h}}}
\newcommand{\vc}{{\vm{c}}}
\newcommand{\viout}{{\aOut{\skew{-2.5}\vm{i}}}}
\newcommand{\viin}{{\aIn{\skew{-2.5}\vm{i}}}}

\newcommand{\vLambda}{{\vm{\lambda}}}

\newcommand{\iu}{{\mkern1.5mu\text{\bfseries i}\mkern1.5mu}}

\begin{document}
\title{A Hybrid Transmission Grid Architecture Enabling Efficient Optimal Power Flow}

\author{Matthias~Hotz,~\IEEEmembership{Student Member,~IEEE,}~and~Wolfgang~Utschick,~\IEEEmembership{Senior Member,~IEEE}%
\thanks{This work was supported by a seed funding of the Munich School of Engineering of Technische Universit\"at M\"unchen within the project TUM.Energy Valley~Bavaria.}%
\thanks{M. Hotz and W. Utschick are with the Department of Electrical and Computer Engineering, Technische Universit\"at M\"unchen, Munich D-80333, Germany (e-mail: matthias.hotz@tum.de, utschick@tum.de).}%
}

\maketitle

\begin{abstract}
The recent rise of electricity generation based on renewable energy sources increases the demand for transmission capacity. Capacity expansion via the upgrade of transmission line capacity, e.g., by conversion to a high-voltage direct current (HVDC) line, is an attractive option. In this paper, it is shown that if the upgrade to HVDC is applied systematically to selected transmission lines across the grid, a hybrid architecture is obtained that enables an efficient and globally optimal solution of the optimal power flow (OPF) problem. More precisely, for conventional meshed AC transmission grids the OPF problem is nonconvex and in general NP-hard, rendering it hard to solve. We prove that after the upgrade to the proposed hybrid architecture, the same mesh topology facilitates an exact convex relaxation of the OPF problem, enabling its globally optimal solution with efficient polynomial time algorithms. This OPF method is then employed in simulations, which demonstrate that the hybrid architecture can increase the effective transmission capacity and substantially reduce the generation costs, even compared to the AC grid with optimal transmission switching.
\end{abstract}


\section{Introduction}
\label{chp:introduction}

\IEEEPARstart{R}{ecently}, many countries experience a significant increase of electricity generation based on renewable energy sources, e.g., wind and solar energy, particularly in Europe~\cite{European-Commission2013a}. This shift in generation leads to an increasingly distributed and fluctuating energy production, which requires the transmission grid to handle and balance strong fluctuations~\cite{Huber2014a, Schaber2012b}. Consequently, the demand for transmission capacity is expected to grow~\cite{Schaber2012b}. As the implementation of new transmission lines is often difficult and protracted, notably due to the obtainment of right of way, the upgrade of transmission line capacity is gaining importance~\cite{Baldick2009a}. To this end, the conversion of AC lines to high-voltage direct current (HVDC) lines is a possible option~\cite{Clerici1991a,Edris2008a,Larruskain2006a,Raju2015a}. HVDC lines feature a controllable power flow and provide dynamic voltage support at the terminals~\cite{Bahrman2007a}, which can further enhance the capacity of the connected AC grid~\cite{Bahrman2007a, Urquidez2012a} and renders it particularly appealing from a technological point of view.

Transmission capacity is not only determined by the available infrastructure, but also by the efficiency of its utilization. The operational task of determining the optimal utilization is known as the \emph{optimal power flow} (OPF) problem~\cite{Frank2012a, Frank2012b}. It identifies the optimal allocation of generation resources (dispatch) and the corresponding state of the transmission grid to serve the load, considering a steady-state AC model of the transmission grid and system constraints~\cite{Wood1996a}. For a general meshed AC transmission grid, this OPF problem is nonconvex and in general NP-hard, thus lacking an efficient computational solution~\cite{Low2014a,Low2014b}. Recently, there has been a lot of research on convex relaxation of the OPF problem to shift it into a computationally tractable domain, cf. the tutorial~\cite{Low2014a, Low2014b} and the references therein. Convex relaxations generally extend the feasible set and, as a consequence, solution recoverability is only ensured under certain conditions, for which the relaxation is called \emph{exact}. It was shown that for radial grids such exact relaxations exist, while for general meshed AC transmission grids no exact relaxations are known yet~\cite{Low2014a, Low2014b, Madani2015a}. However, a mesh topology is vital for the capacity and reliability of a transmission grid, thus these methods are not applicable if recoverability of the solution must be ensured. For this reason, typically simplified system models are employed, e.g., the ``decoupled power flow'' or the widely used ``DC power flow''~\cite{Wood1996a}. Although an OPF formulation based on these simplified models improves computational tractability, the model mismatch requires more conservative system constraints and leads to a suboptimal utilization.

\begin{figure*}[!t]
\centering
\subfloat[]{%
\includegraphics{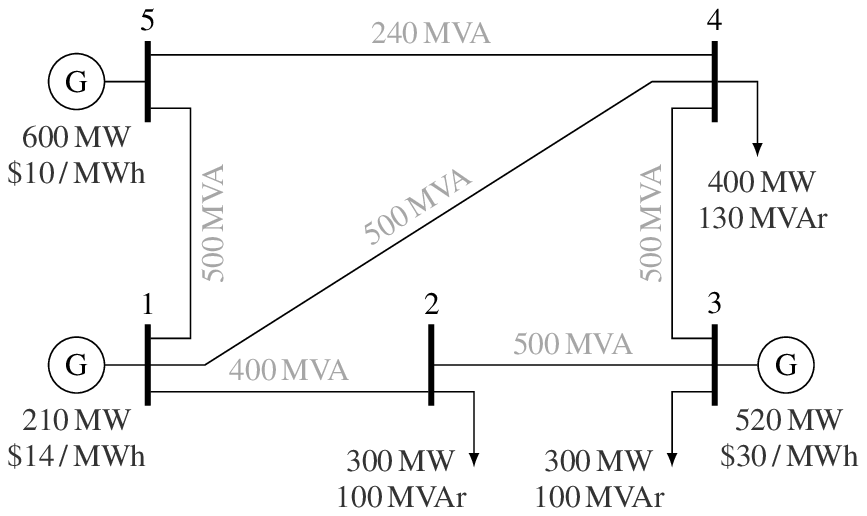}%
\label{fig:sim:grid:acg}}%
\hspace*{6mm}%
\subfloat[]{%
\includegraphics{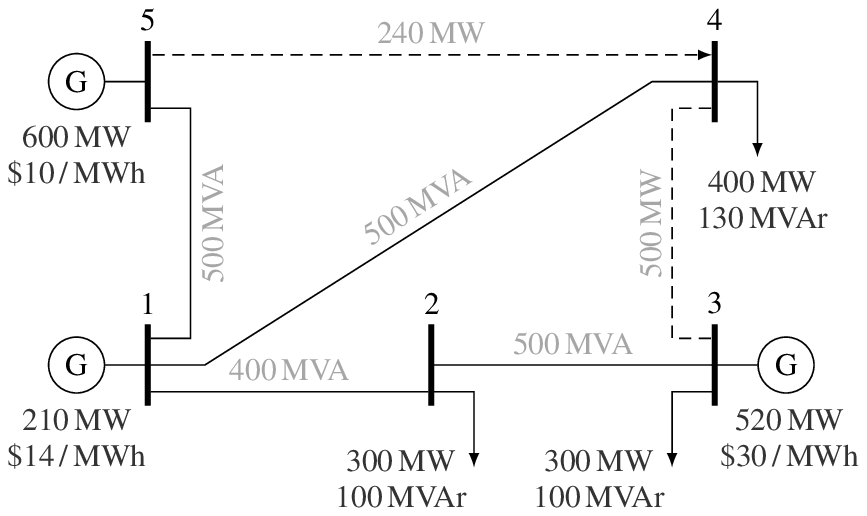}%
\label{fig:sim:grid:htg}}%
\caption{Single-line diagrams of an adapted PJM system (cf.~\cite{Li2010b}), where (a) is the original AC transmission grid and (b) is an exemplary hybrid transmission grid obtained by upgrading AC line $3$\;--\;$4$ and $5$\;--\;$4$ to HVDC. Note that (unidirectional) HVDC lines are represented by a (directional) dashed line.}
\label{fig:sim:grid}
\end{figure*}

Consequently, it is desired to perform the capacity expansion of the transmission infrastructure such that it also supports the grid utilization via an efficient and globally optimal solution of the OPF problem. In this work, a capacity expansion approach for transmission grids is presented, which enables a convex formulation of the OPF problem that is solved globally optimal with efficient polynomial time algorithms. To this end, we propose a hybrid architecture, where an arbitrary spanning tree of the transmission grid is fixed and the AC lines outside this spanning tree are upgraded to HVDC lines as illustrated in Fig.~\ref{fig:sim:grid}. Thus, one AC line in every loop within the transmission grid is upgraded to HVDC. It is important to notice that transmission grids are typically sparsely meshed, i.e., the number of AC lines that cause loops is small compared to the total number of AC lines. Furthermore, as the spanning tree can be chosen arbitrarily, many options are available for the choice of lines that are subject to a capacity upgrade.

The result of the capacity expansion is a hybrid AC/DC grid. The OPF problem for hybrid AC/DC grids was already studied, e.g. in~\cite{Wiget2012a,Feng2014a,Aragues-Penalba2015a,Iggland2015a,Baradar2013a}, but without considering a particular architecture. As it is inherently a nonconvex problem due to the AC power flow equations and global optimization algorithms are computationally intractable for large-scale grids, typically either algorithms that in general converge to local optima~\cite{Wiget2012a,Feng2014a,Aragues-Penalba2015a} or simplified models that lead to a suboptimal utilization~\cite{Iggland2015a,Baradar2013a} are employed. In contrast, the proposed hybrid architecture induces structural properties that enable an exact convex relaxation and, thus, an efficient and globally optimal solution of the OPF problem.

Finally, it shall be pointed out that Farivar and Low~\cite{Farivar2013b} obtained an exact convex relaxation of the OPF problem for meshed AC grids by introducing ideal phase shifting transformers to all lines outside some spanning tree. However, this has the major drawback of requiring potentially large, stability-endangering or even intractable phase shifts at some transformers. In contrast, the proposed hybrid architecture supports stability via the advantages of HVDC technology.

\subsection{Contributions and Outline}

In Section~\ref{chp:model}, the hybrid architecture is introduced via a steady-state system model. To this end, a novel approach is pursued to describe the architecture of this hybrid transmission grid. This allows to preserve its structural properties in the description of the electrical behavior, which is essential for the efficient solution of the OPF problem. Subsequently, in Section~\ref{chp:constraints} an appropriate formulation of system constraints is presented. A particular focus is put on an expressive branch flow constraint to support congestion management. The OPF problem for the hybrid transmission grid and its convex relaxation is established in Section~\ref{chp:opf}, where the exactness of the relaxation is proven in Section~\ref{chp:exactness}. This does not only generalize previous results on convex relaxation of the OPF problem to the meshed hybrid transmission grid but also extends the theory to a more refined system model, which allows the establishment of exactness on simple and intuitive requirements on physical parameters. In Section~\ref{chp:simulation}, the OPF method is demonstrated in simulations, which illustrate that the proposed hybrid architecture can substantially reduce the generation costs and increase the effective transmission capacity. Finally, Section~\ref{chp:conclusion} concludes~the~paper.

\subsection{Notation}
\label{chp:introductionpt1:notation}

The set of real numbers is denoted by $\nR$, the set of nonnegative real numbers by $\nRnn$, the set of complex numbers by $\nC$, and the set of Hermitian matrices in $\nC^{N\times N}$ by $\nS^N$. The imaginary unit is denoted by $\displaystyle\iu=\sqrt{-1}$. For $x\in\nC$, its real part is $\real(x)$, its imaginary part is $\imag(x)$, its absolute value is $\abs{x}$, the principal value of its argument is $\arg(x)\in(-\pi,\pi]$, and its complex conjugate is $x\conj$. For a matrix $\mA$, its transpose is $\mA\tran$, its conjugate (Hermitian) transpose is $\mA\herm$, its trace is $\trace(\mA)$, its rank is $\rank(\mA)$, and the element in row $i$ and column $j$ is $[\mA]_{i,j}$. For two matrices $\mA,\mB\in\nS^N$, $\mA\psd\mB$ denotes that $\mA-\mB$ is positive semidefinite and $\mA\pd\mB$ that $\mA-\mB$ is positive definite. For real-valued vectors, inequalities are component-wise. The vector $\ve_n$ denotes the $n$th standard basis vector of appropriate dimension. For logical values $a$ and $b$, their logical conjunction is $a\;\land\;b$. For a set $\sS$, its cardinality is denoted by $\abs{\sS}$ and its interior by $\interior(\sS)$.

\section{System Model}
\label{chp:model}

The system model for the hybrid transmission grid comprises the description of its \emph{architecture} and the description of its \emph{electrical behavior}. In the following, AC lines, cables, transformers, and phase shifters are referred to as \emph{AC branches}, HVDC lines and cables as well as back-to-back converters are referred to as \emph{DC branches}, and points of interconnection, generation injection, and load connection are called \emph{buses}.

\subsection{Architecture}
\label{chp:model:topology}

The architecture of the hybrid transmission grid is described by the directed multigraph $\gG = {\{\sV, \sE, \sD, \faE, \fbE, \faD, \fbD\}}$, where
\begin{enumerate}
	\item $\sV = \{1,\ldots,\cNV\}$ is the set of buses,
	\item $\sE = \{1,\ldots,\cNE\}$ is the set of AC branches,
	\item $\sD = \{1,\ldots,\cND\}$ is the set of DC branches,
	\item $\faE : \sE\rightarrow\sV$ maps an AC branch to its source bus,
	\item $\fbE : \sE\rightarrow\sV$ maps an AC branch to its destination bus,
	\item $\faD : \sD\rightarrow\sV$ maps a DC branch to its source bus, and
	\item $\fbD : \sD\rightarrow\sV$ maps a DC branch to its destination bus.
\end{enumerate}
The set elements are, without loss of generality, assumed to be consecutive numbers to facilitate the electrical description in matrix notation. The directionality of AC branches is \emph{not} related to the direction of power flow and can be chosen arbitrarily. To characterize the architecture of the hybrid transmission grid, some sets are defined.

\begin{figure}[!t]
\centering
\includegraphics{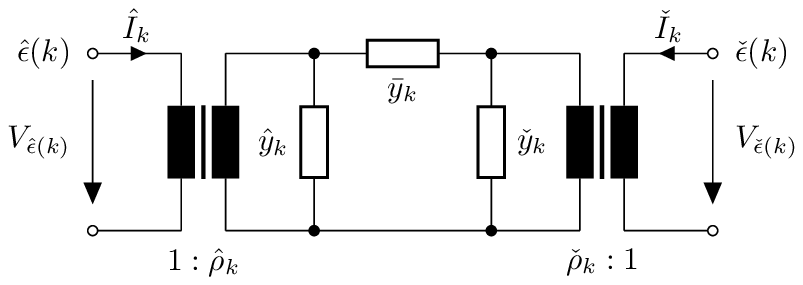}
\vspace{-0.6em}
\caption{Model for AC branch $k\in\sE$, which connects bus $\faE(k)$ and bus $\fbE(k)$.}
\label{fig:acbranch}
\end{figure}

\begin{definitionx}
The set $\aOut{\sBE}(n)\subseteq\sE$ and $\aIn{\sBE}(n)\subseteq\sE$ of AC branches outgoing and incoming at bus $n\in\sV$, respectively,~is
\begin{align}
	\aOut{\sBE}(n) &= \{ k \in \sE : \faE(k) = n \}\allowdisplaybreaks\\
	\aIn{\sBE}(n) &= \{ k \in \sE : \fbE(k) = n \}\;.
\end{align}
\end{definitionx}
\begin{definitionx}
The set $\aOut{\sBD}(n)\subseteq\sD$ and $\aIn{\sBD}(n)\subseteq\sD$ of DC branches outgoing and incoming at bus $n\in\sV$, respectively,~is
\begin{align}
	\aOut{\sBD}(n) &= \{ l \in \sD : \faD(l) = n \}\\
	\aIn{\sBD}(n) &= \{ l \in \sD : \fbD(l) = n \}\;.
\end{align}
\end{definitionx}
\begin{definitionx}
The set $\sEU\subseteq\sV\times\sV$ of undirected edges underlying the AC branches $\sE$ (obtained by removing directionality)~is
\begin{equation}
\sEU = \bigcup_{k\in\sE}\big\{(\faE(k),\fbE(k)),(\fbE(k),\faE(k))\big\}\;.
\end{equation}
\end{definitionx}
In the proposed hybrid architecture, the AC branches must \emph{connect all buses} and form a \emph{tree network topology}. Thus, for a conventional AC grid, this hybrid architecture is obtained by upgrading the AC lines outside some arbitrary spanning tree to HVDC. To keep the model as simple as possible, it is assumed that there do not exist any branches that start and end at the same bus and that there do not exist any parallel AC branches between two buses. The latter may occur in practice, e.g., in the case of parallel transmission lines, but such parallel entities can be modeled as a single AC branch.
\begin{requirement}\label{def:noselfloops}
The multigraph $\gG$ does not comprise any self-loops of length one, i.e., $\nexists\,k\in\sE:\faE(k)=\fbE(k)$ and $\nexists\,l\in\sD:\faD(l)=\fbD(l)$.
\end{requirement}
\begin{requirement}\label{def:noparallelacbranches}
The subgraph $\gG_\text{AC} = {\{\sV, \sE, \faE, \fbE\}}$ with all buses and AC branches does not contain any parallel or antiparallel edges, i.e.,
\begin{align*}
	\text{1)}&\ \ \nexists\,k_1,k_2\in\sE, k_1\neq k_2:\faE(k_1)=\faE(k_2)\;\land\;\fbE(k_1)=\fbE(k_2)\\
	\text{2)}&\ \ \nexists\,k_1,k_2\in\sE, k_1\neq k_2:\faE(k_1)=\fbE(k_2)\;\land\;\fbE(k_1)=\faE(k_2)\;.
\end{align*}
\end{requirement}
\begin{requirement}\label{def:actree}
The undirected graph $\aBar{\gG}_\text{AC} = {\{\sV, \sEU\}}$ underlying the subgraph $\gG_\text{AC} = {\{\sV, \sE, \faE, \fbE\}}$ with all buses and AC branches is a \emph{tree}, i.e., $\aBar{\gG}_\text{AC}$ is \emph{connected} and \emph{acyclic}.
\end{requirement}
\begin{corollary}\label{cor:numacbranches}
There are $\cNE=\cNV-1$ AC branches (cf.~\cite{Diestel2010a}).
\end{corollary}

\subsection{Electrical Model}
\label{chp:model:electrical}

The electrical model for the hybrid transmission grid as presented below is based on \textsc{Matpower}~\cite{Zimmerman2011a} and its mathematical formulation is derived from the work of Low~\cite{Low2014a, Low2014b} and Bose \emph{et~al.}~\cite{Bose2012a, Bose2011a}.

\subsubsection{AC Branch Model}


AC branches are represented via the common branch model in Fig.~\ref{fig:acbranch}. For AC branch $k\in\sE$, it comprises two shunt admittances $\aOut{y_k},\aIn{y_k}\in\nC$, a series admittance $\aBar{y_k}\in\nC$, and two complex voltage ratios $\aOut{\rho_k},\aIn{\rho_k}\in\nC\setminus\{0\}$, where $\abs{\aOut{\rho_k}}$ and $\abs{\aIn{\rho_k}}$ is the tap ratio and $\arg(\aOut{\rho_k})$ and $\arg(\aIn{\rho_k})$ the phase shift of the respective transformer. Let $\vv = [V_1,\ldots,V_{\cNV}]\tran\in\nC^{\cNV}$ be the bus voltage vector, $\viout = [\aOut{I_1},\ldots,\aOut{I_{\cNE}}]\tran\in\nC^{\cNE}$ the source current vector, $\viin = [\aIn{I_1},\ldots,\aIn{I_{\cNE}}]\tran\in\nC^{\cNE}$ the destination current vector, and
\begin{align}
	\aOut{\mY} &= \sum_{k\in\sE}\ve_k(\aOut{\alpha_k}\ve_{\faE(k)}\tran + \aOut{\beta_k}\ve_{\fbE(k)}\tran)\in\nC^{\cNE\times\cNV}\\
	\aIn{\mY} &= \sum_{k\in\sE}\ve_k(\aIn{\alpha_k}\ve_{\fbE(k)}\tran + \aIn{\beta_k}\ve_{\faE(k)}\tran)\in\nC^{\cNE\times\cNV}
\end{align}
the source and destination admittance matrix, in which
\begin{subequations}\label{eqn:alphabeta}
\begin{align}
	\aOut{\alpha_k} &= \abs{\aOut{\rho_k}}^2(\aBar{y_k}+\aOut{y_k}) &
	\aOut{\beta_k}  &= -\rho_k\aBar{y_k}\\
	\aIn{\alpha_k}  &= \abs{\aIn{\rho_k}}^2(\aBar{y_k}+\aIn{y_k}) &
	\aIn{\beta_k}   &= -\rho_k\conj\aBar{y_k}
\end{align}
\end{subequations}
and $\rho_k = \aOut{\rho_k\conj}\aIn{\rho_k}$ is the total voltage ratio. Therewith, the source and destination branch currents are given by
\begin{align}
	\viout = \aOut{\mY}\vv &&
	\viin  = \aIn{\mY}\vv\;.
\end{align}

\subsubsection{Bus Model}

\begin{figure}[!t]
\centering
\includegraphics{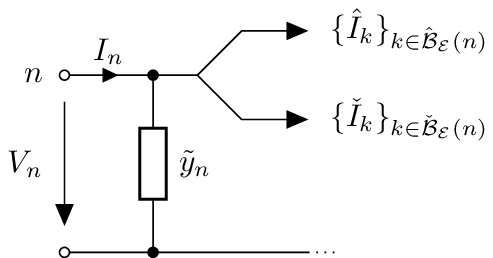}
\caption{Model for bus $n\in\sV$.}
\label{fig:bus}
\vspace{-0.025em}
\end{figure}

Buses are modeled as depicted in Fig.~\ref{fig:bus}. For bus $n\in\sV$, it comprises a shunt admittance $\aShunt{y_n}\in\nC$, connections to the outgoing AC branches $k\in\aOut{\sBE}(n)$ as well as to the incoming AC branches $k\in\aIn{\sBE}(n)$, and an injection port. Let $\vi = [I_1,\ldots,I_{\cNV}]\tran\in\nC^{\cNV}$ be the injection current vector and
\begin{equation}
	\mY = \sum_{n\in\sV}\ve_n\Big[\alpha_n\ve_n\tran\ +\ \sum_{\mathclap{k\in\aOut{\sBE}(n)}}\aOut{\beta_k}\ve_{\fbE(k)}\tran\ +\ \sum_{\mathclap{k\in\aIn{\sBE}(n)}}\aIn{\beta_k}\ve_{\faE(k)}\tran\Big]\in\nC^{\cNV\times\cNV}
\end{equation}
the bus admittance matrix, in which
\begin{equation}
	\alpha_n = \aShunt{y_n}\ +\ \ \sum_{\mathclap{k\in\aOut{\sBE}(n)}}\aOut{\alpha_k}\ +\ \ \sum_{\mathclap{k\in\aIn{\sBE}(n)}}\aIn{\alpha_k}\;.
\end{equation}
Therewith, the injection currents are given by
\begin{equation}
	\vi = \mY\vv\;.
\end{equation}

\subsubsection{DC Branch Model}
\label{chp:model:electrical:dcbranch}

\begin{figure}[!t]
\centering
\includegraphics{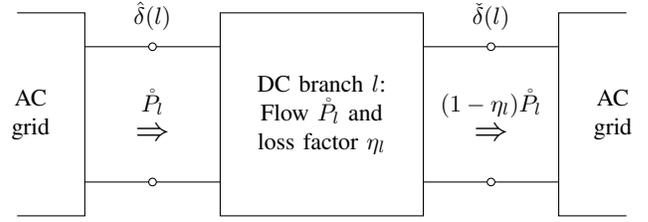}
\caption{Model for DC branch $l\in\sD$, which connects the injection port of the source bus $\faD(l)$ to the injection port of the destination bus $\fbD(l)$.}
\label{fig:dcbranch}
\end{figure}

DC branches are modeled as illustrated in Fig.~\ref{fig:dcbranch}. For DC branch $l\in\sD$, the nonnegative amount $\aDc{P_l}\in\nRnn$ of active power is transferred from source bus (terminal) $\faD(l)$ to destination bus (terminal) $\fbD(l)$ with a loss factor $\eta_l\in[0,1)$. Note that any static losses of the DC branch can be modeled as fixed loads at the adjacent buses and any reactive power injection capabilities, e.g., in case of voltage source converter (VSC) based HVDC lines~\cite{Bahrman2007a}, can be modeled as a (virtual) generator. If required, the AC-side transformers may be modeled explicitly as AC branches. Bidirectionally operable HVDC lines are modeled via antiparallel DC branches.

\subsubsection{Power Injection}

The hybrid transmission grid comprises the combination of DC branches and an AC subgrid. Thus, the total power injection at some bus constitutes the cumulative power injection into both. Let $\vp = \big[\aDc{P_1},\ldots,\aDc{P_{\cND}}\big]\tran\in\nRnn^{\cND}$ be the DC branch flow vector. Therewith, the active power $P_n\in\nR$ and reactive power $Q_n\in\nR$ injected at bus $n\in\sV$ can be quantified as
\begin{align}
	P_n &= \vv\herm\mP_n\vv + \vh_n\tran\vp &\qquad
	Q_n &= \vv\herm\mQ_n\vv
\intertext{in which $\mP_n,\mQ_n\in\nS^{\cNV}$ and $\vh_n\in\nR^{\cND}$ are}
	\mP_n &= (\mS_n+\mS_n\herm)/2 &
	\mQ_n &= (\mS_n-\mS_n\herm)/(2\iu)
		\label{eqn:sysmodel:pnqnmtx}
\end{align}
\begin{equation}
	\vh_n =\ \ \sum_{\mathclap{l\in\aOut{\sBD}(n)}}\ve_l\ -\ \sum_{\mathclap{l\in\aIn{\sBD}(n)}}\,(1-\eta_l)\ve_l
\end{equation}
and $\mS_n = \mY\herm\ve_n\ve_n\tran\in\nC^{\cNV\times\cNV}$.

\subsubsection{Requirements}

The OPF method presented in Section~\ref{chp:opf} is based on certain properties that generally hold in transmission grids. These are established below.
\begin{definition}\label{def:passiveacbranches}
All AC branches are passive, i.e., ${\real(\aBar{y_k})\geq 0}$, ${\real(\aOut{y_k})\geq 0}$, and ${\real(\aIn{y_k})\geq 0}$, for all $k\in\sE$.
\end{definition}
\begin{definition}
The conductance of the series admittance of all AC branches is nonzero, i.e., $\real(\aBar{y_k})\neq 0$, for all $k\in\sE$.
\end{definition}
\begin{corollary}\label{cor:seriesadmittance}
$\real(\aBar{y_k})>0$ and $\abs{\aBar{y_k}}>0$, for all $k\in\sE$.
\end{corollary}
\begin{definition}\label{def:inductiveseriesadmittance}
The susceptance of the series admittance of all AC branches is inductive, i.e., $\imag(\aBar{y_k})\leq 0$, for all $k\in\sE$.
\end{definition}
\begin{definition}\label{def:properinsulation}
All AC branches are properly insulated, i.e., $\abs{\aOut{y_k}}/\abs{\aBar{y_k}}\leq 1$ and $\abs{\aIn{y_k}}/\abs{\aBar{y_k}}\leq 1$, for all $k\in\sE$.
\end{definition}
\begin{definition}\label{def:totalphaseshift}
All AC branches exhibit a total phase shift of less than or equal to $90$\textdegree, i.e., $\abs{\arg(\rho_k)}\leq \pi/2$, for all $k\in\sE$.
\end{definition}

\section{System Constraints}
\label{chp:constraints}

The power flow within the hybrid transmission grid must respect certain constraints as introduced below. Their formulation maintains a quadratic form in the bus voltages, which is essential for the efficient solution of the OPF problem. Furthermore, it departs from the constraint on apparent power flow, which is common to describe the capacity of transmission lines (``MVA rating''). The limit on apparent power flow is actually threefold. For short transmission lines it represents a physical constraint of thermal nature that is proportional to the current, for medium-length transmission lines it casts a stability constraint related to the voltage drop along the line, and for long transmission lines it represents a stability constraint related to the voltage angle difference of the adjacent buses~\cite[Ch.~6.1.12]{Kundur1994a}, \cite[Ch.~4.9]{Bergen2000a}. Here, these three constraints are implemented directly, which does not only increase expressiveness and accuracy of the constraint~\cite{Baldick2009a} but also improves its mathematical structure.

\setcounter{subsubsection}{0}

\subsubsection{Bus Voltage Magnitude}

Transmission grids are designed for a certain voltage range, i.e., at bus $n\in\sV$, the voltage magnitude must satisfy $\abs{V_n}\in[\aLB{V_n},\aUB{V_n}]\subset\nRnn$, with $\aLB{V_n}<\aUB{V_n}$. With $\mM_n = \ve_n\ve_n\tran\in\nS^{\cNV}$, this can be put as
\begin{equation}
	\aLB{V_n}^2 \leq \vv\herm\mM_n\vv \leq \aUB{V_n}^2\;.
\end{equation}

\subsubsection{AC Branch Current Magnitude}

The physical limits on the flow on AC branch $k\in\sE$ are expressed as $\abs{\aOut{I_k}} \leq \aOutUB{I}_k$ and $\abs{\aIn{I_k}} \leq \aInUB{I}_k$, where $\aOutUB{I}_k,\aInUB{I}_k\in\nRnn\setminus\{0\}$. In quadratic form, this renders
\begin{align}
	\vv\herm\mIout_k\vv \leq \aOutUB{I}_k^2 &&
	\vv\herm\mIin_k\vv \leq \aInUB{I}_k^2
\end{align}
where $\mIout_k,\mIin_k\in\nS^{\cNV}$ are
\begin{align}
	\mIout_k=\aOut{\mY}\herm\ve_k\ve_k\tran\aOut{\mY} &&
	\mIin_k=\aIn{\mY}\herm\ve_k\ve_k\tran\aIn{\mY}\;.
\end{align}
\begin{remark}\label{rem:lineflow}
Common line capacity specifications comprise an upper bound on the apparent power flow, i.e.,
\begin{align}
	\abs{\aOut{S}_k} = \abs{\aOut{I}_k}\abs{V_{\faE(k)}} \leq \aUB{S}_k &&\ \quad
	\abs{\aIn{S}_k} = \abs{\aIn{I}_k}\abs{V_{\fbE(k)}} \leq \aUB{S}_k\ \
\end{align}
in which $\aUB{S}_k\in\nRnn\setminus\{0\}$ is the line rating for AC branch $k\in\sE$. If the rating is thermally binding, the bus voltage specification can be utilized to derive the (conservative) substitute bounds
\begin{align}
	\abs{\aOut{I_k}} \leq \aOutUB{I}_k = \aUB{S}_k / \aUB{V_{\faE(k)}} &&\qquad
	\abs{\aIn{I_k}}  \leq \aInUB{I}_k  = \aUB{S}_k / \aUB{V_{\fbE(k)}}\;.\
\end{align}
\end{remark}

\subsubsection{DC Branch Power Flow}

The power flow on a DC branch is limited by its physical capabilities and may need to maintain some minimum flow for proper operation. With $\aLB{\vp} = \big[\aDcLB{P}_1,\ldots,\aDcLB{P}_{\cND}\big]\tran,\ \aUB{\vp} = \big[\aDcUB{P}_1,\ldots,\aDcUB{P}_{\cND}\big]\tran\in\nRnn^{\cND}$, this can be expressed collectively for all DC branches as
\begin{equation}
	\aLB{\vp} \leq \vp \leq \aUB{\vp}\;.
\end{equation}

\subsubsection{AC Branch Voltage Magnitude Drop}

The stability-related limit on the relative voltage magnitude drop $\nu_k\in\nR$ along AC branch $k\in\sE$, i.e.,
\begin{equation}
	\nu_k = \abs{V_{\fbE(k)}}/\abs{V_{\faE(k)}} - 1
\end{equation}
reads $\nu_k\in[\aLB{\nu_k},\aUB{\nu_k}]\subset[-1,\infty)$, with $\aLB{\nu_k}<\aUB{\nu_k}$, or equivalently
\begin{align}
	\vv\herm\aLB{\mM_k}\vv \leq 0 &&
	\vv\herm\aUB{\mM_k}\vv \leq 0
\end{align}
in which $\aLB{\mM_k},\aUB{\mM_k}\in\nS^{\cNV}$ are
\begin{align}
	\aLB{\mM_k} &= (1 + \aLB{\nu_k})^2\mM_{\faE(k)} - \mM_{\fbE(k)}\\
	\aUB{\mM_k} &= \mM_{\fbE(k)} - (1 + \aUB{\nu_k})^2\mM_{\faE(k)}\;.
\end{align}

\subsubsection{AC Branch Voltage Angle Difference}

The stability-related limit on the voltage angle difference $\delta_k\in\nR$ along AC branch $k\in\sE$, i.e.,
\begin{equation}
	\delta_k = \arg(V_{\faE(k)}\conj V_{\fbE(k)})
\end{equation}
reads $\delta_k\in[\aLB{\delta_k},\aUB{\delta_k}]\subset(-\pi/2,\pi/2)$, with $\aLB{\delta_k}<\aUB{\delta_k}$. Note that
\begin{align}
	\real(V_{\faE(k)}\conj &V_{\fbE(k)}) \geq 0 \\
	\hspace*{-0.2em}
	\tan(\aLB{\delta_k}) \leq \imag(V_{\faE(k)}\conj V_{\fbE(k)})\,/\,&\real(V_{\faE(k)}\conj V_{\fbE(k)}) \leq \tan(\aUB{\delta_k})
\end{align}
is an equivalent formulation of this constraint, which can be written in quadratic form as
\begin{align}
	\vv\herm\mA_k\vv \leq 0
	&&
	\vv\herm\aLB{\mA_k}\vv \leq 0
	&&
	\vv\herm\mAub_k\vv \leq 0
\end{align}
where $\mA_k,\aLB{\mA_k},\mAub_k\in\nS^{\cNV}$ are $\mA_k = -\aOut{\mM_k} - \aOut{\mM_k\herm}$ and
\begin{align}
	\aLB{\mA_k} &= (\tan(\aLB{\delta_k}) + \iu)\aOut{\mM_k} + (\tan(\aLB{\delta_k}) - \iu)\aOut{\mM_k\herm}\\
	\mAub_k &= -(\tan(\aUB{\delta_k}) + \iu)\aOut{\mM_k} - (\tan(\aUB{\delta_k}) - \iu)\aOut{\mM_k\herm}
\end{align}
with $\aOut{\mM_k} = \ve_{\faE(k)}\ve_{\fbE(k)}\tran\in\nS^{\cNV}$. To efficiently solve the OPF problem, $\aLB{\delta_k} $ and $\aUB{\delta_k}$ are required to satisfy
\begin{align}
	\hspace*{-0.55em}
	-\pi/2 < \aLB{\delta_k} \leq -\arg(\rho_k)
	\ \,\text{and}\ \,
	-\arg(\rho_k) \leq \aUB{\delta_k} < \pi/2\;.
	\label{eqn:sysmodel:deltabndcrt}
\end{align}
As the total phase shift $\arg(\rho_k)$ is usually smaller in magnitude than common values of $\aLB{\delta_k} $ and $\aUB{\delta_k}$, which are about $\pm 40$\textdegree\ to $\pm 50$\textdegree~(cf. e.g.~\cite{Bergen2000a}), this generally holds in practice.

\subsubsection{Bus Power Injection}

At bus $n\in\sV$, let the load exhibit an active power demand $\aPload{P_n}\in\nR$ and reactive power demand $\aPload{Q_n}\in\nR$ and let the connected generation utility provide an active power injection range $[\aPgenmin{P_n},\aPgenmax{P_n}]\subset\nR$ and reactive power injection range $[\aPgenmin{Q_n},\aPgenmax{Q_n}]\subset\nR$. Then, the net injection range amounts to $P_n\in[\aLB{P_n},\aUB{P_n}]$ and $Q_n\in[\aLB{Q_n},\aUB{Q_n}]$ with $\aLB{P_n} = \aPgenmin{P_n} - \aPload{P_n}$, $\aLB{Q_n} = \aPgenmin{Q_n} - \aPload{Q_n}$, $\aUB{P_n} = \aPgenmax{P_n} - \aPload{P_n}$, $\aUB{Q_n} = \aPgenmax{Q_n} - \aPload{Q_n}$. Thus, the power injection constraint reads
\begin{align}
	\aLB{P_n} \leq \vv\herm\mP_n\vv + \vh_n\tran\vp \leq \aUB{P_n} \\
	\aLB{Q_n} \leq \vv\herm\mQ_n\vv \leq \aUB{Q_n}\;.\ \,
\end{align}

\section{Optimal Power Flow}
\label{chp:opf}

The optimal power flow in the hybrid transmission grid is considered with respect to the cost- and/or loss-minimizing allocation of generation facilities, where the resulting power flows are compliant with the system constraints. It is assumed that a proper unit commitment~\cite{Wood1996a, Padhy2004a} was performed a priori, i.e., the generators' on/off status is preassigned.

\subsection{Objective}
\label{chp:opf:objective}

Let the generated active power $\aPgen{P_n}=P_n+\aPload{P_n}$ at bus $n\in\sV$ be associated with a cost $\gamma_n\in\nRnn$ in dollar per watt-hour (\$/Wh). Then, the hourly generation cost is given by
\begin{equation}
	\aCost{f}(\vv,\vp) = \vv\herm\aCost{\mC}\vv + \aCost{\vc}\tran\vp + \Gamma_\text{L}
\end{equation}
where $\aCost{\mC}\in\nS^{\cNV}$, $\aCost{\vc}\in\nR^{\cND}$, and $\Gamma_\text{L}\in\nR$ are
\begin{align}
	\hspace*{-0.6em}
	\aCost{\mC} = \sum_{n\in\sV} \gamma_n\mP_n\ \;\quad
	\aCost{\vc} = \sum_{n\in\sV} \gamma_n\vh_n\ \;\quad
	\Gamma_\text{L} = \sum_{n\in\sV} \gamma_n\aPload{P_n}.
	\label{eqn:costobjparamdef}
\end{align}
The electrical loss comprises the losses on AC branches, DC branches, and at buses. These amount to
\begin{equation}
	\aLoss{f}(\vv,\vp) = \vv\herm\aLoss{\mC}\vv + \aLoss{\vc}\tran\vp
\end{equation}
where $\aLoss{\mC}\in\nS^{\cNV}$ and $\aLoss{\vc}\in\nR^{\cND}$ are
\begin{align}
	\aLoss{\mC} = \sum_{k\in\sE}\aLoss{\mP_k} + \sum_{n\in\sV}\real(\aShunt{y_n})\mM_n &&\qquad
	\aLoss{\vc} = \sum_{l\in\sD}\eta_l\ve_l
	\label{eqn:lossobjparamdef}
\end{align}
with $\aLoss{\mP_k} = (\aLoss{\mS_k}+\aLoss{\mS_k\herm})/2\in\nS^{\cNV}$ and $\aLoss{\mS_k} = \aOut{\mY}\herm\ve_k\ve_{\faE(k)}\tran+\aIn{\mY}\herm\ve_k\ve_{\fbE(k)}\tran\in\nC^{\cNV\times\cNV}$. The objective of the OPF problem is to minimize the weighted sum of generation cost and electrical loss, i.e.,
\begin{equation}
	f(\vv,\vp)
		= w\aCost{f}(\vv,\vp) + \aLoss{\gamma}\aLoss{f}(\vv,\vp)
		= \vv\herm\mC\vv + \vc\tran\vp + w\Gamma_\text{L}\ 
	\label{eqn:objfunction}
\end{equation}
in which $\mC\in\nS^{\cNV}$ and $\vc\in\nR^{\cND}$ are
\begin{align}
	\mC = w\aCost{\mC} + \aLoss{\gamma}\aLoss{\mC} &&\qquad
	\vc = w\aCost{\vc} + \aLoss{\gamma}\aLoss{\vc}\quad
\end{align}
and $w,\aLoss{\gamma}\in\nRnn$ are the weights for the generation cost and electrical loss, respectively. To obtain an objective function for a representative (virtual) total hourly cost, the weight $w$ is considered dimensionless and $\aLoss{\gamma}$ is in dollar per watt-hour (\$/Wh), i.e., an (artificial) cost for the electrical loss. In order to ensure an efficient solution of the OPF problem, the weight for the loss term is required to be nonzero (cf. Appendix~\ref{apx:thm:objmtxelements}). Note that this is not limiting in~practice, as the ratio of weights is not restricted.
\begin{definition}\label{def:lossweightnonzero}
The weight for the loss term is strictly positive, i.e., $\aLoss{\gamma}>0$.
\end{definition}

\subsection{Optimal Power Flow Problem}

The OPF problem to determine the optimal utilization of the hybrid transmission grid is the minimization of the above objective while considering the system constraints,~i.e.,
\begin{subequations}\label{eqn:opfoptprobfull}
\begin{align}
&\minimize_{\mathclap{\vv\in\nC^{\cNV},\,\vp\in\nR^{\cND}}} & \vv\herm\mC\vv &+ \vc\tran\vp\\
&\subjectto
&\vv\herm\mM_n\vv &\leq \aUB{V_n}^2\,, & \forall n\in\sV
	\label{eqn:opfoptprobfull:crt:mn}\\
&&\vv\herm\left(-\mM_n\right)\vv &\leq -\aLB{V_n}^2\,, & \forall n\in\sV\\
&&\vv\herm\mP_n\vv + \vh_n\tran\vp &\leq \aUB{P_n}\,, & \forall n\in\sV
	\label{eqn:opfoptprobfull:crt:pn}\\
&&\vv\herm\mQ_n\vv &\leq \aUB{Q_n}\,, & \forall n\in\sV\\
&&\vv\herm\mIout_k\vv &\leq \aOutUB{I}_k^2\,, & \forall k\in\sE\\
&&\vv\herm\mIin_k\vv &\leq \aInUB{I}_k^2\,, & \forall k\in\sE\\
&&\vv\herm\aLB{\mM_k}\vv &\leq 0\,, & \forall k\in\sE\\
&&\vv\herm\aUB{\mM_k}\vv &\leq 0\,, & \forall k\in\sE\\
&&\vv\herm\mA_k\vv &\leq 0\,, & \forall k\in\sE\\
&&\vv\herm\aLB{\mA_k}\vv &\leq 0\,, & \forall k\in\sE\\
&&\vv\herm\mAub_k\vv &\leq 0\,, & \forall k\in\sE
	\label{eqn:opfoptprobfull:crt:akub}\\
&&\vp &\leq \aUB{\vp}\,,
	\label{eqn:opfoptprobfull:crt:pub}\\
&&-\vp &\leq -\aLB{\vp}\;.
	\label{eqn:opfoptprobfull:crt:plb}
\end{align}
\end{subequations}
This optimization problem identifies the state of the hybrid transmission grid, i.e., the bus voltage vector $\vv$ and DC branch flow vector $\vp$, that minimizes the weighted sum of generation cost and electrical loss. The corresponding allocation of generation capacity at the individual buses follows immediately from the respective power injection. It should be pointed out that the lower bounds on active and reactive power injection are omitted, which has technical reasons as discussed later on.

The OPF problem~\eqref{eqn:opfoptprobfull} can be cast as a quadratically constrained quadratic program (QCQP) in $\vv$ and $\vp$, and due to the fact that $-\mM_n$ is negative semidefinite and $\mP_n$, $\mQ_n$, $\aLB{\mM_k}$, $\aUB{\mM_k}$, $\mA_k$, $\aLB{\mA_k}$, and $\mAub_k$ are in general indefinite, this optimization problem is nonconvex and NP-hard to solve (cf.~\cite{Boyd2004a, Zheng2011a, Low2014a}). To circumvent this issue and shift the OPF problem into a computationally tractable domain, it is reformulated as a convex optimization problem below, which requires that~\eqref{eqn:opfoptprobfull} is essentially strictly feasible.
\begin{assumption}\label{ass:strictfeas}
There exists a $\vv\in\nC^{\cNV}$ and $\vp\in\nRnn^{\cND}$ so that~\eqref{eqn:opfoptprobfull:crt:mn} to~\eqref{eqn:opfoptprobfull:crt:akub} hold with strict inequality and~\eqref{eqn:opfoptprobfull:crt:pub} and~\eqref{eqn:opfoptprobfull:crt:plb} holds (componentwise) either strict or with equality.
\end{assumption}

It should be noted that the solution of the OPF problem~\eqref{eqn:opfoptprobfull} is not unique. In particular, if $(\aOpt{\vv},\aOpt{\vp})$ is an optimizer of~\eqref{eqn:opfoptprobfull}, then $(\aShunt{\vv},\aOpt{\vp})$ with $\aShunt{\vv}={\aOpt{\vv}\en^{\iu\phi}}$ and $\phi\in\nR$ is also an optimizer. This is usually resolved by designating some bus $n_\text{ref}\in\sV$ as the reference bus and adding the constraint $\arg(\ve_{n_\text{ref}}\tran\vv) = \phi_\text{ref}$. However, in the following this constraint becomes redundant and, hence, is omitted in~\eqref{eqn:opfoptprobfull}.

\subsection{Convex Relaxation}

Nonconvex QCQPs arise in many engineering problems and semidefinite relaxation (SDR) is an established technique to efficiently solve them or approximate their solution~\cite{Luo2010a}. To apply SDR to the OPF problem of the hybrid transmission grid, the notation is simplified by rewriting~\eqref{eqn:opfoptprobfull} as
\begin{subequations}\label{eqn:opfoptprobabst}
\begin{align}
	&\minimize_{\mathclap{\vv\in\nC^{\cNV},\,\vp\in\nR^{\cND}}\hspace{-0.5em}} && \vv\herm\mC_0\vv + \vc_0\tran\vp\\
	&\subjectto
	&& \vv\herm\mC_m\vv + \vc_m\tran\vp \leq b_m\,,\ m = 1,...\,,M
		\label{eqn:opfoptprobabst:crt}
\end{align}
\end{subequations}
which comprises $M=4\cNV+7\cNE+2\cND$ constraints with a corresponding parametrization. A key step in SDR is the utilization of the cyclic property of the trace~\cite{Luo2010a}, i.e.,
\begin{equation*}
	\vv\herm\mC_m\vv = \trace(\vv\herm\mC_m\vv) = \trace(\mC_m\vv\vv\herm) = \trace(\mC_m\mV)
\end{equation*}
in which $\mV={\vv\vv\herm}$. Therewith,~\eqref{eqn:opfoptprobabst} can be reformulated in $\mV$, where $\mV$ must be Hermitian, positive semidefinite, and have rank~$1$ to facilitate the decomposition into ${\vv\vv\herm}$. Due to the linearity of the trace, the objective as well as the inequality constraints are linear in $\mV$. As the set of positive semidefinite matrices is convex, the nonconvexity of the problem stems solely from the rank constraint. In SDR, the optimization problem is rendered \emph{convex} by dropping the rank constraint~\cite{Luo2010a}, i.e., the \emph{relaxed} OPF problem is obtained~as
\begin{subequations}\label{eqn:opfrelaxed}
\begin{align}
	&\minimize_{\mathclap{\mV\in\nS^{\cNV},\,\vp\in\nR^{\cND}}\hspace{-0.5em}} && \!\trace(\mC_0\mV) + \vc_0\tran\vp\\
	&\subjectto
	&&\!\trace(\mC_m\mV) + \vc_m\tran\vp \leq b_m\,,\ m = 1,...\,,M
		\label{eqn:opfrelaxed:crtineq}\\
	&&&\mV\psd\mZero\;.
		\label{eqn:opfrelaxed:crtpsd}
\end{align}
\end{subequations}
Note that the relaxation eliminates the ambiguity in the bus voltage angle, as 
$[\mV]_{i,j}$ comprises $V_i V_j\conj$ and, thus, involves only angle differences. It is also important to recognize that the nonconvex problem~\eqref{eqn:opfoptprobfull} has $\cNV+\cND$ variables, whereas~\eqref{eqn:opfrelaxed} has $\cNV^2+\cND$ variables. However, the matrices $\mC_m$, $m=0,\ldots,M$, are chordal sparse and only $\cNV+2\cNE+\cND={3\cNV+\cND-2}$ variables in~\eqref{eqn:opfrelaxed} are actually multiplied by nonzero coefficients (cf. Appendix~\ref{apx:thm:crtmtxcone} and Corollary~\ref{cor:numacbranches}). This chordal sparsity pattern can be utilized to maintain efficient solvability for large-scale grids, cf. e.g.~\cite{Klerk2010a, Fukuda2001a, Nakata2003a, Burer2003a, Andersen2010a}.

\subsection{Solution Recovery}

If \emph{all} optimizers of~\eqref{eqn:opfrelaxed} have rank~1 the relaxation is called \emph{exact}.\footnote{Note that the definition of exactness varies, e.g., in~\cite{Low2014b} it is more stringent.} Then, a solution $(\aOpt{\vv},\aOpt{\vp})$ of~\eqref{eqn:opfoptprobfull} can be recovered from a solution $(\aOpt{\mV},\aOpt{\vp})$ of~\eqref{eqn:opfrelaxed} by decomposing $\aOpt{\mV}$ into an outer product $\aOpt{\mV}={\aOpt{\vv}(\aOpt{\vv})\herm}$, e.g., via the method in~\cite[Sec.~IV-D]{Low2014a} or an eigenvalue decomposition (EVD). For the latter, denote the eigenvalues of $\aOpt{\mV}$ as $\sigma_1\geq\sigma_2\geq\ldots\geq\sigma_{\cNV}$ and the corresponding eigenvectors as $\aOpt{\vv}_1,\ldots,\aOpt{\vv}_{\cNV}$. Then, $\aOpt{\vv}=\sqrt{\sigma_1}\aOpt{\vv}_1$. Note that in case of exactness and a numerically precise solution, $\sigma_2=\ldots=\sigma_{\cNV}=0$. However, due to the limited accuracy of practical solvers these are typically nonzero. Thus, it is reasonable to define a reconstruction error measure
\begin{equation}
	\kappa(\aOpt{\mV}) = \abs{\sigma_2/\sigma_1}
\end{equation}
which is zero in case of exact reconstruction and increases with increasing reconstruction error. After the solution recovery, the absolute angle at the reference bus $n_\text{ref}$ can be matched to the voltage angle reference $\phi_\text{ref}$ via
\begin{equation}
	\aShunt{\vv} = \aOpt{\vv}\exp(\iu[\phi_\text{ref}-\arg(\ve_{n_\text{ref}}\tran\aOpt{\vv})])\;.
\end{equation}

\section{Exactness of the Relaxation}
\label{chp:exactness}

In the following, it is proven that the convex relaxation~\eqref{eqn:opfrelaxed} of the OPF problem~\eqref{eqn:opfoptprobfull} is \emph{exact} for \emph{all} hybrid transmission grids that comply with the system model in Section~\ref{chp:model}. Thus, the operational task of determining the optimal utilization is performed efficiently by solving the relaxed OPF problem~\eqref{eqn:opfrelaxed} using polynomial time algorithms.

In summary, it is proven that the Karush-Kuhn-Tucker (KKT) conditions of~\eqref{eqn:opfrelaxed} with Slater's constraint qualification (cf. e.g.~\cite{Boyd2004a,Bazaraa2006a}), which are necessary for optimality, only permit optimizers with rank~1. It should be pointed out that this proof was inspired by the discussion of exact convex relaxation of QCQPs with underlying acyclic graph structure by Bose \emph{et al.}~\cite{Bose2012a},\footnote{See also the abridged version~\cite{Bose2014a} and their earlier result~\cite{Bose2011a}.} which is based on the influential work of Lavaei and Low in~\cite{Lavaei2012a}. They combine a result of van der Holst~\cite[Th.~3.4]{Holst2003a} with the vanishing inner product of the optimization variable and some matrix to prove that the optimizer has rank~1 in case of a positive semidefinite objective matrix and some further preconditions. This remarkable connection is used here as well, however, embedded in a very different proof technique. In contrast to the dual problem based approach in~\cite{Bose2012a}, the proof below relies on the KKT conditions of the relaxed OPF problem. Furthermore, a novel method to ensure a tree graph behind the matrix emerging in the inner product is introduced, which enables the utilization of the unique properties of the hybrid transmission grid and does \emph{not} require a positive semidefinite objective matrix.

\subsection{Mathematical Preliminaries}

A major part of the proof is based on cones and half-spaces in the complex plane, where the required mathematical framework and some properties are established below.
\begin{definitionx}
The \emph{conic hull} $\cone(\sS)\subseteq\nC$ of a set $\sS\subset\nC$ with finite cardinality $\abs{\sS}$ is defined as
\begin{align*}
	\cone(\sS) = \bigg\{ x\in\nC : x = &\sum_{i=1}^{j}\alpha_i x_i,\ x_i\in\sS,\\[-0.25em]
		&\quad\alpha_i\in\nRnn,\ i,j\in\{1,\ldots,\abs{\sS}\}\bigg\}\;.
\end{align*}
\end{definitionx}
\begin{corollary}\label{cor:convhullarg}
The conic hull of any nonempty set $\sS\subset\nC$, where $\sS\neq\{0\}$, can be described via the arguments of the nonzero elements in $\sS$, in particular,
\begin{equation*}
	\cone(\sS) = \cone\big(\{\en^{\iu\arg(x)}\}_{x\in\sS:\abs{x}>0}\big)\;.
\end{equation*}
\end{corollary}
\begin{definitionx}
The \emph{complex conjugate} $\sS\conj\subseteq\nC$ of a set $\sS\subseteq\nC$ is defined as $\sS\conj = \{ x\in\nC : x\conj\in\sS \}$\,.
\end{definitionx}
\begin{definitionx}
The set $\sS\subset\nC$ is a \emph{half-space} with \emph{normal} $p\in\nC\setminus\{0\}$ if $\sS = \{x\in\nC : \real(p\conj x)\leq 0\}$\,.
\end{definitionx}
\begin{definition}
The \emph{interior} $\interior(\sS)$ of a half-space $\sS\subset\nC$ with normal $p\in\nC\setminus\{0\}$ is $\interior(\sS) = \{x\in\nC : {\real(p\conj x) < 0}\}$\,.
\end{definition}
\begin{corollary}\label{cor:halfplaneoriginprop}
If $\sS\subset\nC$ is a half-space, then $0\notin\interior(\sS)$.
\end{corollary}
\begin{proposition}\label{prp:scalednormal}
Let $\sS\subset\nC$ be a half-space with 
normal $p\in\nC\setminus\{0\}$. Then, $\aBar{p}=\alpha p$, with $\alpha>0$, is also a normal of~$\sS$.
\end{proposition}
\begin{IEEEproof}
$\real(\aBar{p}\conj x)=\alpha\real(p\conj x)\leq 0\,,\ \forall x\in\sS\,$.
\end{IEEEproof}
\begin{proposition}
Let $\sS\subset\nC$ be a half-space with normal $p\in\nC\setminus\{0\}$. Then, $\sS\conj$ is a half-space and $p\conj$ a normal of $\sS\conj$.
\end{proposition}
\begin{IEEEproof}
Let $\aBar{x}\in\sS\conj$ be the complex conjugate of $x\in\sS$, i.e., $\aBar{x}=x\conj$. For all $\aBar{x}\in\sS\conj$ and $\aBar{p}=p\conj$,
\begin{equation*}
	\real(\aBar{p}\conj\aBar{x})=\real(p x\conj)=\real((p\conj x)\conj)=\real(p\conj x)\leq 0\;.\IEEEQEDhereeqn
\end{equation*}
\end{IEEEproof}
\begin{lemma}\label{lem:halfplanesuminint}
Let $\sS\subset\nC$ be a half-space, $x_1\in\sS$, $x_2\in\interior(\sS)$, and $\aBar{x}=x_1+x_2$. Then, $\aBar{x}\in\interior(\sS)$.
\end{lemma}
\begin{IEEEproof}
Let $p\in\nC\setminus\{0\}$ be a normal of $\sS$. By definition, $\real(p\conj x_1)\leq 0$ and $\real(p\conj x_2)<0$. Therefore,
\begin{equation*}
	\real(p\conj\aBar{x})
		=\real(p\conj x_1)+\real(p\conj x_2)<0\;.\IEEEQEDhereeqn
\end{equation*}
\end{IEEEproof}

\subsection{Sufficient Condition for Rank-1 Solutions}
\label{chp:exactness:kkt}

The proof is separated into two parts, the establishment of a sufficient condition for rank~1 solutions, as covered in this section, and the fulfillment of this condition, which follows in Section~\ref{chp:exactness:fulfillment}. To start with, the KKT conditions and Slater's constraint qualification for the relaxed OPF problem~\eqref{eqn:opfrelaxed} read as follows (cf. e.g.~\cite[Ch.~5.9]{Boyd2004a}). Let ${\fPsiCone:\nRnn^M\rightarrow\nS^{\cNV}}$ and ${\fpsiCone:\nRnn^M\rightarrow\nR^{\cND}}$ be defined as
\begin{align}
	&\fPsiCone(\vLambda)=\sum_{m=1}^M[\vLambda]_m\mC_m &\qquad
	\fpsiCone(\vLambda)=\sum_{m=1}^M[\vLambda]_m\vc_m
\intertext{and $\fPsi:\nRnn^M\rightarrow\nS^{\cNV}$ and $\fpsi:\nRnn^M\rightarrow\nR^{\cND}$ be defined as}
	&\fPsi(\vLambda)=\mC_0 + \fPsiCone(\vLambda) &\qquad
	\fpsi(\vLambda)=\vc_0 + \fpsiCone(\vLambda)\;.
		\label{eqn:kktpsidef}
\end{align}
The KKT conditions for~\eqref{eqn:opfrelaxed} are given by the primal feasibility in~\eqref{eqn:opfrelaxed:crtineq} and~\eqref{eqn:opfrelaxed:crtpsd}, the dual feasibility
\begin{align}
	\fPsi(\vLambda)\psd\mZero\label{eqn:KKT:DF:Psi} &&
	\fpsi(\vLambda)=\mZero
\end{align}
and the complementary slackness condition
\begin{align}
	[\vLambda]_m\big(\trace(\mC_m\mV) + \vc_m\tran\vp - b_m\big)=\;&0\,,\ \ m = 1,\ldots,M
		\label{eqn:KKT:CS:ineqcrt}\\
	\trace(\fPsi(\vLambda)\mV)=\;&0\label{eqn:KKT:CS:trace}
\end{align}
where $\vLambda\in\nRnn^M$ is the Lagrangian multiplier. Slater's constraint qualification is fulfilled if there exists a feasible tuple $(\mV,\vp)$ with $\mV\pd\mZero$, which renders the KKT conditions necessary for optimality. This is established by the following theorem that is proven in Appendix~\ref{apx:thm:strictfeas}.
\begin{theorem}\label{thm:strictfeas}
Consider the convex relaxation~\eqref{eqn:opfrelaxed} of the optimization problem~\eqref{eqn:opfoptprobfull}. If Assumption~\ref{ass:strictfeas} holds, there exists a feasible tuple $(\mV,\vp)$ in~\eqref{eqn:opfrelaxed} for which $\mV\pd\mZero$.
\end{theorem}

Thus, the KKT conditions are necessary for optimality in~\eqref{eqn:opfrelaxed} and it is shown in the following that due to the complementary slackness condition in~\eqref{eqn:KKT:CS:trace} the rank of $\fPsi(\vLambda)$ restricts the rank of $\mV$. To this end, consider the following lemma, which follows from the rank-nullity theorem.
\begin{lemma}\label{lem:basisandev}
Let $\mA,\mB\in\nS^N$ be positive semidefinite matrices. Then, $\trace(\mA\mB)=0$ implies $\rank(\mA)+\rank(\mB)\leq N$.
\end{lemma}

It follows from~\eqref{eqn:opfrelaxed:crtpsd}, \eqref{eqn:KKT:DF:Psi}, \eqref{eqn:KKT:CS:trace}, and Lemma~\ref{lem:basisandev} that
\begin{equation}
	\rank(\mV) \leq \cNV - \rank(\fPsi(\vLambda))\;.
\end{equation}
Thus, if $\rank(\fPsi(\vLambda)) \geq \cNV-1$ then any potential optimizer of~\eqref{eqn:opfrelaxed} has rank $1$ and, consequently, the relaxation is exact. Note that ${\rank(\mV)}=0$ implies that all bus voltages are zero, which means a shutdown of the grid and, thus, is not feasible if the OPF problem is properly specified. To state a condition on $\fPsi(\vLambda)$ that ensures its rank is greater than or equal to $\cNV-1$, a result of van der Holst is utilized (cf.~\cite{Bose2012a}), which requires the notion of an \emph{associated} graph.
\begin{definitionx}[cf.~\cite{Holst2003a}]\label{def:mtxgraph}
Let $\mA\in\nS^N$ be a Hermitian matrix. The undirected graph $\fMtxGraph{\mA}=(\sV_{\mA},\sE_{\mA})$ associated with $\mA$ comprises the vertices $\sV_{\mA}=\{1,\ldots,N\}$ and the edges
\begin{equation*}
	\sE_{\mA}=\{(i,j)\in\sV_{\mA}\times\sV_{\mA}:[\mA]_{i,j}\neq 0\;\land\;i\neq j\,\}\;.
\end{equation*}
\end{definitionx}

The result~\cite[Th.~3.4]{Holst2003a} of van der Holst states that $\rank(\fPsi(\vLambda)) \geq {\cNV-1}$ if $\fPsi(\vLambda)\psd\mZero$ and $\fMtxGraph{\fPsi(\vLambda)}$ is a \emph{tree}. Due to~\eqref{eqn:KKT:DF:Psi}, $\fPsi(\vLambda)$ is positive semidefinite for all KKT points and, consequently, a \emph{sufficient condition} for the exactness of the relaxation is that $\fMtxGraph{\fPsi(\vLambda)}$ is a tree, for all $\vLambda\in\nRnn^M$.

\subsection{Fulfillment of the Rank-1 Condition}
\label{chp:exactness:fulfillment}

In the following, it is shown that, for all $\vLambda\in\nRnn^M$, the graph $\fMtxGraph{\fPsi(\vLambda)}$ is equivalent to $\aBar{\gG}_\text{AC}$, i.e., $\fMtxGraph{\fPsi(\vLambda)}\equiv\aBar{\gG}_\text{AC}$. Then, as $\aBar{\gG}_\text{AC}$ is a tree by construction (see Definition~\ref{def:actree}), it follows that $\fMtxGraph{\fPsi(\vLambda)}$ is a tree, thus the sufficient condition for rank~1 solutions is fulfilled and the relaxation is exact. To prove $\fMtxGraph{\fPsi(\vLambda)}\equiv\aBar{\gG}_\text{AC}$, it is shown that $\fMtxGraph{\fPsi(\vLambda)}$ inherits all edges of $\aBar{\gG}_\text{AC}$ by showing that those and only those corresponding off-diagonal elements in $\fPsi(\vLambda)$ are nonzero, i.e., the nonzero elements of $\fPsi(\vLambda)$ coincide with nonzero elements of the adjacency matrix of $\aBar{\gG}_\text{AC}$ (off the diagonal). More precisely, ensuring that the nonzero elements in $\fPsi(\vLambda)$ are a subset of those in the adjacency matrix of $\aBar{\gG}_\text{AC}$ ensures acyclicity, whereas showing that they do not vanish ensures connectivity.
\begin{definitionx}
The cone $\sC_k$ associated with AC branch ${k\in\sE}$ is defined as $\sC_k = \cone(\sX_k)$, where
\begin{align}
\begin{split}
	\sX_k = \{
		\aOut{\beta_{k}}/2,\,
		&\aIn{\beta_{k}}\conj/2,
		-\aOut{\beta_{k}}/(2\iu),\,
		\aIn{\beta_{k}}\conj/(2\iu),\,
		\aOut{\alpha_{k}\conj}\aOut{\beta_{k}},\,
		\\&
		\aIn{\alpha_{k}}\aIn{\beta_{k}\conj},\,
		-1,\,
		\tan(\aLB{\delta_k}) + \iu,\,
		-\tan(\aUB{\delta_k}) - \iu \}\,.
\end{split}
\end{align}
\end{definitionx}
These cones allow the characterization of potentially nonzero off-diagonal elements in $\fPsiCone(\vLambda)$, where all other off-diagonal elements are zero, i.e., those not associated with any AC~branch.
\begin{theorem}\label{thm:crtmtxcone}
For all $\vLambda\in\nRnn^M$, the off-diagonal elements of $\fPsiCone(\vLambda)$, i.e., $[\fPsiCone(\vLambda)]_{i,j}$ for $i,j\in\sV$ and $i\neq j$, satisfy
\begin{align*}
	[\fPsiCone(\vLambda)]_{i,j}
		&\in\left\{
			\begin{array}{ll}
				\sC_{k_1}
					& \text{if}\ \exists\,k_1\in\sE:\faE(k_1)=i\;\land\;\fbE(k_1)=j\\
				\sC_{k_2}\conj
					& \text{if}\ \exists\,k_2\in\sE:\faE(k_2)=j\;\land\;\fbE(k_2)=i\\
				\{0\} & \text{otherwise}\;.
			\end{array}
		\right.
\end{align*}
\end{theorem}
\begin{IEEEproof}
See Appendix~\ref{apx:thm:crtmtxcone}.
\end{IEEEproof}
This characterization can be linked to half-spaces in $\nC$.
\begin{theorem}\label{thm:crtconeishalfplane}
For all $k\in\sE$, the set $\sC_k$ is a subset of a half-space $\sH_k\subset\nC$ with the normal $p=\rho_k$, i.e.,
\begin{equation}
	\sC_k\subseteq\sH_k = \{x\in\nC : \real(\rho_k\conj x)\leq 0\}\;.
	\label{eqn:hndef}
\end{equation}
\end{theorem}
\begin{IEEEproof}
See Appendix~\ref{apx:thm:crtconeishalfplane}.
\end{IEEEproof}
The nonzero elements of $\mC_0$ are inside these same half-spaces.
\begin{theorem}\label{thm:objmtxelements}
For all $\vLambda\in\nRnn^M$, the off-diagonal elements of $\mC_0$, i.e., $[\mC_0]_{i,j}$ for $i,j\in\sV$ and $i\neq j$, satisfy
\begin{align*}
	[\mC_0]_{i,j}
		&\in\left\{
			\begin{array}{ll}
				\interior(\sH_{k_1})
					& \text{if}\ \exists\,k_1\in\sE:\faE(k_1)=i\;\land\;\fbE(k_1)=j\\
				\interior(\sH_{k_2}\conj)
					& \text{if}\ \exists\,k_2\in\sE:\faE(k_2)=j\;\land\;\fbE(k_2)=i\\
				\{0\} & \text{otherwise}\;.
			\end{array}
		\right.
\end{align*}
\end{theorem}
\begin{IEEEproof}
See Appendix~\ref{apx:thm:objmtxelements}.
\end{IEEEproof}
This allows to characterize the off-diagonal elements of $\fPsi(\vLambda)$.
\begin{theorem}\label{thm:psimtxelements}
For all $\vLambda\in\nRnn^M$, the off-diagonal elements of $\fPsi(\vLambda)$, i.e., $[\fPsi(\vLambda)]_{i,j}$ for $i,j\in\sV$ and $i\neq j$, satisfy
\begin{align*}
	[\fPsi(\vLambda)]_{i,j}
		&\in\left\{
			\begin{array}{ll}
				\nC\setminus\{0\}
					& \text{if}\ \exists\,k\in\sE:\faE(k),\fbE(k)\in\{i,j\}\\
				\{0\} & \text{otherwise}\;.
			\end{array}
		\right.
\end{align*}
\end{theorem}
%
\begin{table}[!t]
\renewcommand{\arraystretch}{\tabstretch}
\caption{Bus-related Data for $n\in\sV=\{1,2,3,4,5\}$}
\label{tab:mdl:bus}
\centering
\begin{tabular}{|c||c|c||c|c|S[table-format=4]|S[table-format=4]||c|}
\hline
\multirow{2}{*}{$n$} & $\aPload{P_n}$ & $\aPload{Q_n}$ & $\aPgenmin{P_n}$ & $\aPgenmax{P_n}$ & {$\aPgenmin{Q_n}$} & {$\aPgenmax{Q_n}$} & $\gamma_n$ \\ \cline{2-8}
& MW & MVAr & MW & MW & {MVAr} & {MVAr} & \$/MWh \\ \hline\hline
1 &  -- &  -- &  0 & 210 & -155 &  155 & 14 \\
2 & 300 & 100 & -- &  -- & {--} & {--} & -- \\
3 & 300 & 100 &  0 & 520 & -390 &  390 & 30 \\
4 & 400 & 130 & -- &  -- & {--} & {--} & -- \\
5 &  -- &  -- &  0 & 600 & -450 &  450 & 10 \\ \hline
\end{tabular}
\end{table}
%
\begin{IEEEproof}
Due to~\eqref{eqn:kktpsidef}, $[\fPsi(\vLambda)]_{i,j}=[\mC_0]_{i,j}+[\fPsiCone(\vLambda)]_{i,j}$ and consider $i\neq j$. It follows from Theorem~\ref{thm:crtmtxcone} and~\ref{thm:objmtxelements} that if there does not exist any AC branch between bus $i$ and $j$, i.e., $\nexists\,k\in\sE:\faE(k),\fbE(k)\in\{i,j\}$, then $[\fPsiCone(\vLambda)]_{i,j}=0$, $[\mC_0]_{i,j}=0$, and thus $[\fPsi(\vLambda)]_{i,j}=0$. If there exists an AC branch $k\in\sE$ between bus $i$ and $j$, i.e., $\exists\,k\in\sE:\faE(k),\fbE(k)\in\{i,j\}$, then it follows from Theorem~\ref{thm:crtmtxcone},~\ref{thm:crtconeishalfplane}, and~\ref{thm:objmtxelements} that $[\fPsiCone(\vLambda)]_{i,j}\in\sS$ and $[\mC_0]_{i,j}\in\interior(\sS)$, where $\sS\subset\nC$ is a half-space, i.e., $\sH_k$ or $\sH_k\conj$. Therewith, Lemma~\ref{lem:halfplanesuminint} states $[\fPsi(\vLambda)]_{i,j}\in\interior(\sS)$ and Corollary~\ref{cor:halfplaneoriginprop} that $[\fPsi(\vLambda)]_{i,j}\neq 0$.
\end{IEEEproof}

Theorem~\ref{thm:psimtxelements} states that $[\fPsi(\vLambda)]_{i,j}$ (with $i\neq j$) is nonzero if there exists an AC branch between bus $i$ and $j$ and zero otherwise. Applying Definition~\ref{def:mtxgraph} yields $\fMtxGraph{\fPsi(\vLambda)} = {\{\sV, \sEU\}}$, which leads to the following corollary and completes the proof of exactness for the relaxation~\eqref{eqn:opfrelaxed} of the OPF problem~\eqref{eqn:opfoptprobfull}.
\begin{corollary}
$\fMtxGraph{\fPsi(\vLambda)}\equiv\aBar{\gG}_\text{AC}$, for all $\vLambda\in\nRnn^M$.
\end{corollary}

\section{Simulation Results}
\label{chp:simulation}

\begin{table}[!t]
\renewcommand{\arraystretch}{\tabstretch}
\caption{AC Branch-related Data for $k\in\sE=\{1,2,3,4\}$}
\label{tab:mdl:acb}
\centering
\begin{tabular}{|c||c|c||c|c||c|}
\hline
\multirow{2}{*}{$k$} & \multirow{2}{*}{$\faE(k)$} & \multirow{2}{*}{$\fbE(k)$} & $\aBar{z_k}=1/\aBar{y_k}$ & $\aOut{y_k},\aIn{y_k}$ & $\aUB{S}_k$ \\ \cline{4-6}
& & & p.u. & p.u. & MVA \\ \hline\hline
1 & 1 & 2 & $0.00281 + \iu 0.02810$ & $\iu 0.00356$ & 400 \\
2 & 1 & 4 & $0.00304 + \iu 0.03040$ & $\iu 0.00329$ & 500 \\
3 & 1 & 5 & $0.00064 + \iu 0.00640$ & $\iu 0.01563$ & 500 \\
4 & 2 & 3 & $0.00108 + \iu 0.01080$ & $\iu 0.00926$ & 500 \\ \hline
\noalign{\vskip 0.5em}
\multicolumn{6}{l}{\hspace*{-0.5em}Converted AC lines:} \\[0.15em] \hline
-- & 3 & 4 & $0.00297 + \iu 0.02970$ & $\iu 0.00337$ & 500 \\
-- & 5 & 4 & $0.00297 + \iu 0.02970$ & $\iu 0.00337$ & 240 \\ \hline
\end{tabular}
\end{table}

\begin{table}[!t]
\renewcommand{\arraystretch}{\tabstretch}
\caption{DC Branch-related Data for $l\in\sD=\{1,2,3\}$}
\label{tab:mdl:dcb}
\centering
\begin{tabular}{|c||c|c||S[table-format=1.1]||c|c|}
\hline
\multirow{2}{*}{$l$} & \multirow{2}{*}{$\faD(l)$} & \multirow{2}{*}{$\fbD(l)$} & $\eta_l$ & $\aDcLB{P}_l$ & $\aDcUB{P}_l$ \\ \cline{4-6}
& & & \parbox{7mm}{\centering \%} & \parbox{7mm}{\centering MW} & \parbox{7mm}{\centering MW} \\ \hline\hline
1 & 3 & 4 & 3.5 & 0 & 500 \\
2 & 4 & 3 & 3.5 & 0 & 500 \\
3 & 5 & 4 & 3.5 & 0 & 240 \\ \hline
\end{tabular}
\end{table}

In the following, the upgrade to the hybrid architecture is discussed and, using the proposed OPF method, the resulting hybrid transmission grid is compared to the AC grid in terms of economic efficiency and effective transmission capacity. To this end, the PJM 5-bus system in~\cite{Li2010b} is considered, where the generators at bus~1 (``Alta'' and ''Park City'') are unified, the generator at bus~4 (``Sundance'') is disconnected to intensify the strain on the grid, and the unconstrained lines are equipped with a rating of $500$\,MVA to accentuate congestion, cf. Fig.~\ref{fig:sim:grid:acg}. This adapted PJM system, where bus~3 serves as the reference bus, is utilized as the reference AC grid. To compare its economic efficiency to the hybrid transmission grid, three different OPF methods with increasing computational complexity are considered for the reference AC grid. Typically, the OPF of AC grids is based on simplified system models to reduce the computational effort. For example, the widely used ``DC power flow'' approximation facilitates an OPF formulation called ``DC OPF''~\cite{Zimmerman2011a} that can be solved in polynomial time, where the required slack power to compensate for the model mismatch is determined by a subsequent AC power flow. This OPF method (``DC\;OPF~\&~AC\;PF'') is used as a reference for the comparison of generation costs. As a second method, the nonconvex OPF based on the AC power flow is considered (``AC\;OPF''), which is accurate but NP-hard to solve. Finally, as a third method, optimal transmission switching~\cite{Fisher2008a} based on the AC\;OPF is applied (``AC\;OTS''), which improves the economic efficiency of an AC grid by switching off selected lines. AC\;OTS, as considered in this work, enumerates all line switching configurations that maintain a connected AC grid, performs an AC\;OPF for every configuration, and selects the one that exhibits the minimum generation cost. Thus, it illustrates the performance of the reference AC grid for the optimal topology, but involves a substantial computational effort as it augments the NP-hard AC\;OPF with a combinatorial problem.

\subsection{Upgrade to the Hybrid Architecture}

\begin{table}[!t]
\renewcommand{\arraystretch}{\tabstretch}
\caption{Case Study 1: OPF Result for the Hybrid Transmission Grid}
\label{tab:res:case1:htg}
\vspace{-\baselineskip}
\begin{minipage}[t]{0.348\textwidth}\vspace{0pt}
\hspace*{0.25mm}
\begin{tabular}{|c||S[table-format=3.2]|S[table-format=3.2]||S[table-format=2.3]|S[table-format=2.3]|}
\hline
{\multirow{2}{*}{$n$}} & {$\aPgen{P_n}$} & {$\aPgen{Q_n}$} & {$\abs{V_n}$} & {$\arg(V_n)$} \\ \cline{2-5}
& {MW} & {MVAr} & {p.u.} & {deg.} \\ \hline\hline
1 & 210.00 & 155.00 & 1.096 &  5.929 \\
2 &   0.00 &   0.00 & 1.088 &  0.555 \\
3 & 201.99 & 227.89 & 1.100 &  0.000 \\
4 &   0.00 &   0.00 & 1.046 &  1.466 \\
5 & 600.00 &  22.94 & 1.100 &  7.441 \\ \hline
\end{tabular}
\end{minipage}
\begin{minipage}[t]{0.135\textwidth}\vspace{0pt}
\begin{tabular}{|c||S[table-format=3.2]|}
\hline
{\multirow{2}{*}{$l$}} & {$\aDc{P_l}$} \\ \cline{2-2}
& {MW} \\ \hline\hline
1 &   0.00 \\
2 &   3.61 \\
3 & 100.66 \\ \hline
\end{tabular}
\end{minipage}\\[0.9\baselineskip]
\hspace*{0.25mm}
\begin{tabular}{|c||S[table-format=3.2]|S[table-format=2.3]|S[table-format=3.2]|S[table-format=2.3]||S[table-format=2.2]|S[table-format=2.3]|}
\hline
{\multirow{2}{*}{$k$}} & {$\abs{\aOut{S}_k}$} & {\multirow{2}{*}{$\aOut{\cPF}_k$}} & {$\abs{\aIn{S}_k}$} & {\multirow{2}{*}{$\aIn{\cPF}_k$}} & {$\nu_k$} & {$\delta_k$} \\ \cline{2-2}\cline{4-4}\cline{6-7}
& {MVA} & & {MVA} & & {\%} & {deg.} \\ \hline\hline
1 & 398.56 &  1.000 & 395.73 & -0.998 & -0.70 & -5.374 \\
2 & 348.45 &  0.888 & 332.90 & -0.921 & -4.55 & -4.463 \\
3 & 498.20 & -1.000 & 499.86 &  0.999 &  0.36 &  1.511 \\
4 & 159.11 &  0.596 & 159.03 & -0.594 &  1.07 & -0.555 \\ \hline
\end{tabular}
\end{table}

Kirchhoff's matrix tree theorem (cf. e.g.~\cite{Merris1994a}) states that the reference AC grid in Fig.~\ref{fig:sim:grid:acg} comprises $11$ spanning trees, i.e., eleven options are available for the upgrade to the hybrid architecture. To determine the upgrade that facilitates the most economic operation, the OPF problem is solved for every option. In the OPF formulation, the AC lines outside the spanning tree are replaced by HVDC lines with an exemplary loss factor of $3.5\%$ and a capacity that coincides with the AC line rating. The uprating due to the conversion is omitted intentionally to focus on the influence of architecture. From this analysis, it follows that the hybrid architecture with minimum generation cost comprises the upgrade of AC line $3$\;--\;$4$ and $5$\;--\;$4$, cf. Fig.~\ref{fig:sim:grid:htg}. The data for this hybrid transmission grid is documented in Table~\ref{tab:mdl:bus}, \ref{tab:mdl:acb}, and~\ref{tab:mdl:dcb}, with $\aLB{V_n}=0.9$\;p.u., $\aUB{V_n}=1.1$\;p.u., $\aShunt{y_n}=0$\;p.u., for all $n\in\sV=\{1,2,3,4,5\}$, and $\aOut{\rho_k}=1$, $\aIn{\rho_k}=1$, $\aUB{\nu_k}=5\%$, $\aLB{\nu_k}=-5\%$, $\aUB{\delta_k}=50$\textdegree, $\aLB{\delta_k}=-50$\textdegree, for all $k\in\sE=\{1,2,3,4\}$. The bounds $\aOutUB{I}_k$ and $\aInUB{I}_k$ are derived from $\aUB{S}_k$ as described in Remark~\ref{rem:lineflow}. For illustration purposes, the HVDC line $5$\;--\;$4$ is considered directional and, for the sake of simplicity, both HVDC lines are assumed to have no static losses and no reactive power injection capabilities.

\begin{table}[!t]
\renewcommand{\arraystretch}{\tabstretch}
\caption{Case Study 1: Dispatch for the Reference AC Grid}
\label{tab:res:case1:acg}
\centering
\begin{tabular}{|c||S[table-format=3.2]|S[table-format=4.2]||S[table-format=3.2]|S[table-format=4.2]||S[table-format=3.2]|S[table-format=4.2]|}
\hline
& \multicolumn{2}{c||}{DC\;OPF \& AC\;PF} & \multicolumn{2}{c||}{AC\;OPF} & \multicolumn{2}{c|}{AC\;OTS} \\ \cline{2-7}
$n$ & $\aPgen{P_n}$ & $\aPgen{Q_n}$ & $\aPgen{P_n}$ & $\aPgen{Q_n}$ & $\aPgen{P_n}$ & $\aPgen{Q_n}$ \\ \cline{2-7}
& {MW} & {MVAr} & {MW} & {MVAr} & {MW} & {MVAr} \\ \hline\hline
1 & 210.00 &  91.97 & 210.00 &  155.00 & 210.00 &  155.00 \\
3 & 328.76 & 259.07 & 324.71 &  378.50 & 274.28 &  390.00 \\
5 & 466.51 &  24.02 & 470.37 & -161.66 & 523.41 & -145.75 \\ \hline
\end{tabular}
\vspace*{-0.25em}
\end{table}

\begin{table}[!t]
\renewcommand{\arraystretch}{\tabstretch}
\caption{Case Study 1: Comparison of Total Generation Cost}
\label{tab:res:case1:cost}
\centering
\begin{tabular}{|l||c|c|c|c|}
\hline
& \multicolumn{3}{c|}{Reference AC Grid} & \multirow{3}{15mm}{\centering Hybrid {\hspace*{-0.15em}Transmission} Grid} \\ \cline{2-4}
& \multirow{2}{13.1mm}{\centering DC\;OPF \& AC\;PF} & \multirow{2}{13.1mm}{\centering AC\;OPF} & \multirow{2}{13.1mm}{\centering AC\;OTS} & \\[-0.4em]
& & & & \\ \hline\hline
Cost      & \costnum{17468} & \costnum{17385} & \costnum{16403} & \costnum{15000} \\ \hline
Reduction &        $0.00\%$ &        $0.47\%$ &        $6.10\%$ &       $14.13\%$ \\ \hline
\end{tabular}
\end{table}

\subsection{Case Study 1: Economic Efficiency}
	\label{chp:simulation:routing}

To identify the optimal utilization of the hybrid transmission grid, the relaxed OPF problem~\eqref{eqn:opfrelaxed} is solved using CVX~\cite{Grant2015a}, where the weights of the objective are set to $w=1$ and $\aLoss{\gamma}=10^{-6}$\;\$/MWh to minimize generation cost. The results are shown in Table~\ref{tab:res:case1:htg} $(\kappa(\aOpt{\mV}) = 4.61\cdot 10^{-12})$, in which the power flow on AC branches is illustrated via the magnitude of the source and destination apparent power flow $\aOut{S}_k = \aOut{I}_k\conj V_{\faE(k)}$ and ${\aIn{S}_k = \aIn{I}_k\conj V_{\fbE(k)}}$ with the corresponding power factors $\aOut{\cPF}_k = \real(\aOut{S}_k)/\abs{\aOut{S}_k}$ and $\aIn{\cPF}_k = \real(\aIn{S}_k)/\abs{\aIn{S}_k}$. For the reference AC grid, \textsc{Matpower}\footnote{Default settings with solver ``OT'' for DC\;OPF and ``MIPS'' for AC\;OPF.}~\cite{Zimmerman2015a} was utilized to implement the three OPF methods. The resulting dispatch is documented in Table~\ref{tab:res:case1:acg}. In case of the DC\;OPF \& AC\;PF, the model mismatch causes a violation of the rating of AC line $5$\;--\;$4$ by $4.11$\,MVA, which shows the necessity of more conservative system constraints for this OPF method.

The total generation cost for the hybrid transmission grid and the reference AC grid is shown in Table~\ref{tab:res:case1:cost}. It can be observed that the reference AC grid exhibits the worst economic efficiency for DC\;OPF \& AC\;PF. AC\;OPF slightly improves the efficiency, whereas AC\;OTS achieves a significant cost reduction by switching off AC line $3$\;--\;$4$, but at a substantial computational cost. In contrast, the hybrid transmission grid achieves more than twice this cost reduction, even though the proposed OPF method is in the same complexity class as the DC\;OPF. It further indicates that the hybrid architecture can mitigate the economic impact of congestion, which is supported by the fact that higher loss factors of the HVDC lines do not significantly affect this result, e.g., for $\eta_l=7.0\%$ the cost reduction is $13.55\%$.

\subsection{Case Study 2: Transmission Capacity}
	\label{chp:simulation:capacity}

To investigate the effective transmission capacity, the active power demand at bus~2 and~4 is increased by $150$\,MW, i.e., $\aPload{P_2}$ and $\aPload{P_4}$ is set to $450$\,MW and $550$\,MW, respectively. The OPF method for the hybrid transmission grid provides the results in Table~\ref{tab:res:b} $(\kappa(\aOpt{\mV}) = 7.33\cdot 10^{-11})$, which corresponds to a total generation cost of \costnum{24044}. In contrast, all three OPF methods for the reference AC grid fail to find a solution. Although this does not imply infeasibility, it still suggests that the reference AC grid cannot serve these loads. Considering that no uprating was performed when the AC lines $3$\;--\;$4$ and $5$\;--\;$4$ were converted to HVDC, this result indicates that the hybrid architecture can increase the effective transmission capacity via its sophisticated power flow control.

\begin{table}[!t]
\renewcommand{\arraystretch}{\tabstretch}
\caption{Case Study 2: OPF Result for the Hybrid Transmission Grid}
\label{tab:res:b}
\vspace{-\baselineskip}
\begin{minipage}[t]{0.348\textwidth}\vspace{0pt}
\hspace*{0.25mm}
\begin{tabular}{|c||S[table-format=3.2]|S[table-format=3.2]||S[table-format=2.3]|S[table-format=2.3]|}
\hline
{\multirow{2}{*}{$n$}} & {$\aPgen{P_n}$} & {$\aPgen{Q_n}$} & {$\abs{V_n}$} & {$\arg(V_n)$} \\ \cline{2-5}
& {MW} & {MVAr} & {p.u.} & {deg.} \\ \hline\hline
1 & 210.00 & 155.00 & 1.095 &  2.345 \\
2 &   0.00 &   0.00 & 1.088 & -0.998 \\
3 & 503.48 & 204.61 & 1.100 &  0.000 \\
4 &   0.00 &   0.00 & 1.041 & -3.413 \\
5 & 600.00 &  39.24 & 1.100 &  3.651 \\ \hline
\end{tabular}
\end{minipage}
\begin{minipage}[t]{0.135\textwidth}\vspace{0pt}
\begin{tabular}{|c||S[table-format=3.2]|}
\hline
{\multirow{2}{*}{$l$}} & {$\aDc{P_l}$} \\ \cline{2-2}
& {MW} \\ \hline\hline
1 &   0.00 \\
2 &   0.00 \\
3 & 166.70 \\ \hline
\end{tabular}
\end{minipage}\\[0.9\baselineskip]
\hspace*{0.25mm}
\begin{tabular}{|c||S[table-format=3.2]|S[table-format=2.3]|S[table-format=3.2]|S[table-format=2.3]||S[table-format=2.2]|S[table-format=2.3]|}
\hline
{\multirow{2}{*}{$k$}} & {$\abs{\aOut{S}_k}$} & {\multirow{2}{*}{$\aOut{\cPF}_k$}} & {$\abs{\aIn{S}_k}$} & {\multirow{2}{*}{$\aIn{\cPF}_k$}} & {$\nu_k$} & {$\delta_k$} \\ \cline{2-2}\cline{4-4}\cline{6-7}
& {MVA} & & {MVA} & & {\%} & {deg.} \\ \hline\hline
1 & 248.71 &  0.999 & 247.00 & -1.000 & -0.69 & -3.343 \\
2 & 431.59 &  0.913 & 410.27 & -0.948 & -5.00 & -5.758 \\
3 & 433.55 & -0.997 & 435.07 &  0.996 &  0.42 &  1.306 \\
4 & 227.25 & -0.893 & 228.79 &  0.889 &  1.12 &  0.998 \\ \hline
\end{tabular}
\end{table}

\section{Conclusion}
\label{chp:conclusion}

In this paper, a hybrid transmission grid architecture was presented that enables the efficient solution of the OPF problem. This hybrid architecture is established by a particular capacity expansion approach, where the transmission lines that are subject to a capacity upgrade via conversion to HVDC are selected such that one AC line in every loop is upgraded. A detailed system model for the hybrid transmission grid as well as an appropriate OPF formulation was introduced. It was proven that the proposed hybrid architecture enables an exact convex relaxation of the OPF problem, where the globally optimal solution can be determined with efficient polynomial time algorithms. Finally, the application of the OPF method was illustrated for an exemplary hybrid transmission grid, where the results show that this hybrid architecture can enable a more economic operation and increase the effective capacity.

This concept also poses further research questions. While the capacity expansion comprises the establishment of the hybrid architecture, the optimal grid-specific and load variation robust choice of lines that are subject to the upgrade remains open. Furthermore, such a capacity expansion constitutes a long-term investment in the infrastructure, necessitating an analysis of realizability and profitability that considers the present grid topology, the cost and properties of the utilized HVDC technology, different load profiles, and so forth.

\appendices

\section{Proof of Theorem~\ref{thm:strictfeas}}
	\label{apx:thm:strictfeas}

Assumption~\ref{ass:strictfeas} states that there exists a feasible tuple $(\vv,\vp)$ in~\eqref{eqn:opfoptprobfull}, where, for all $m = 1,\ldots,M$, the constraints satisfy
\begin{equation}
	\vv\herm\mC_m\vv + \vc_m\tran\vp
		\ \left\{
			\begin{array}{ll}
				< b_m & \text{if}\ \mC_m\neq\mZero\\
				\leq b_m & \text{if}\ \mC_m=\mZero\;.
			\end{array}
		\right.
\end{equation}
Consider the tuple $(\mV,\vp)$ with $\mV = \vv\vv\herm + \varepsilon\mI\,$, where $\varepsilon>0$ is an arbitrary positive scalar and $\mI$ is the $\cNV\times\cNV$ identity matrix, which implies $\mV\pd\mZero$. It follows that
\begin{align}
	\trace(\mC_m\mV) + \vc_m\tran\vp
		&= \trace(\mC_m(\vv\vv\herm + \varepsilon\mI)) + \vc_m\tran\vp\\
		&= \vv\herm\mC_m\vv + \vc_m\tran\vp + \varepsilon\trace(\mC_m)\;.
\end{align}
As $\trace(\mC_m)=0$ if $\mC_m=\mZero$, there exists some $\varepsilon>0$ such that $\trace(\mC_m\mV) + \vc_m\tran\vp \leq b_m\,$, for all $m=1,\ldots,M$. Therefore, $(\mV,\vp)$ is feasible in~\eqref{eqn:opfrelaxed}, which completes the proof.

\section{Proof of Theorem~\ref{thm:crtmtxcone}}
	\label{apx:thm:crtmtxcone}

\begin{table*}[!t]
\renewcommand{\arraystretch}{\tabstretch}
\setlength{\tabcolsep}{1.5em}
\caption{Elements of the Constraint Matrices for AC Branch $k\in\sE$}
\label{tab:acbmtxelements}
\centering
\begin{tabular}{|c|c||c|c|c|c|c|c|c|}
\hline
$i$ & $j$ & $[\mIout_k]_{i,j}$ & $[\mIin_k]_{i,j}$ & $[\aLB{\mM_k}]_{i,j}$ & $[\aUB{\mM_k}]_{i,j}$ & $[\mA_k]_{i,j}$ & $[\aLB{\mA_k}]_{i,j}$ & $[\mAub_k]_{i,j}$ \\
\hline\hline
$\faE(k)$ & $\faE(k)$ & $\abs{\aOut{\alpha_k}}^2$            & $\abs{\aIn{\beta_k}}^2$            & $(1 + \aLB{\nu_k})^2$       & $-(1 + \aUB{\nu_k})^2$ & $0$  &                          $0$ &                          $0$ \\
$\faE(k)$ & $\fbE(k)$ & $\aOut{\alpha_k\conj}\aOut{\beta_k}$ & $\aIn{\alpha_k}\aIn{\beta_k\conj}$ & $0$                         & $0$                    & $-1$ & $\tan(\aLB{\delta_k}) + \iu$ & $-\tan(\aUB{\delta_k}) - \iu$ \\
$\fbE(k)$ & $\faE(k)$ & $\aOut{\alpha_k}\aOut{\beta_k\conj}$ & $\aIn{\alpha_k\conj}\aIn{\beta_k}$ & $0$                         & $0$                    & $-1$ & $\tan(\aLB{\delta_k}) - \iu$ & $-\tan(\aUB{\delta_k}) + \iu$\\
$\fbE(k)$ & $\fbE(k)$ & $\abs{\aOut{\beta_k}}^2$             & $\abs{\aIn{\alpha_k}}^2$           & $-1$                        & $1$                    & $0$  &                          $0$ &                          $0$ \\ \hline
\multicolumn{2}{|c||}{otherwise} & 0 & 0 & 0 & 0 & 0 & 0 & 0 \\ \hline
\end{tabular}
\end{table*}

As $\vLambda\in\nRnn^M$, the element $[\fPsiCone(\vLambda)]_{i,j}$ of $\fPsiCone(\vLambda)$ is a \emph{conic combination} of the elements $[\mC_m]_{i,j}$ of the matrices $\mC_m$, i.e.,
\begin{equation}
	[\fPsiCone(\vLambda)]_{i,j}=\sum_{m=1}^M[\vLambda]_m[\mC_m]_{i,j}
		\label{eqn:kktpsiconeelements}
\end{equation}
for all $i,j\in\sV$. Consequently, $[\fPsiCone(\vLambda)]_{i,j}$ is in the conic hull of the elements $[\mC_m]_{i,j}$, i.e.,
\begin{equation}
	[\fPsiCone(\vLambda)]_{i,j}\in\cone(\{[\mC_m]_{i,j}\}_{m=1,\ldots,M})\;.
\end{equation}
Due to the parametrization for the OPF problem, the elements $[\mC_m]_{i,j}$ amount to the elements of the matrices $\mM_n$, $-\mM_n$, $\mP_n$, $\mQ_n$, $\mIout_k$, $\mIin_k$, $\aLB{\mM_k}$, $\aUB{\mM_k}$, $\mA_k$, $\aLB{\mA_k}$, and $\mAub_k$, for all $n\in\sV$ and $k\in\sE$. By construction, $\mM_n$ is diagonal and for $\mP_n$ and $\mQ_n$ it follows from~\eqref{eqn:sysmodel:pnqnmtx} that their elements are given by $[\mP_n]_{i,j} = ([\mS_n]_{i,j}+[\mS_n]\conj_{j,i})/2$ and $[\mQ_n]_{i,j} = ([\mS_n]_{i,j}-[\mS_n]\conj_{j,i})/(2\iu)$, in which
\begin{align*}
	[\mS_n]_{i,j}
		= \ve_i\tran\mS_n\ve_j
		= \ve_i\tran\mY\herm\ve_n\ve_n\tran\ve_j
		= \left\{
			\begin{array}{ll}
				[\mY]_{j,i}\conj & \text{if}\ j=n\\
				0 & \text{otherwise}
			\end{array}
		\right.
\end{align*}
and $[\mY]_{j,i}\conj = \ve_i\tran\mY\herm\ve_j$ is determined as\footnote{Note that due to Definition~\ref{def:noselfloops} and~\ref{def:noparallelacbranches} the nonzero off-diagonal elements are related to a unique AC branch.}
\begin{align*}
	[\mY]_{j,i}\conj
		&= \sum_{n\in\sV}\ve_i\tran\Big[\alpha_n\conj\ve_n +\ \sum_{\mathclap{k\in\aOut{\sBE}(n)}}\aOut{\beta_k\conj}\ve_{\fbE(k)} +\ \sum_{\mathclap{k\in\aIn{\sBE}(n)}}\aIn{\beta_k\conj}\ve_{\faE(k)}\Big]\ve_n\tran\ve_j\allowdisplaybreaks\\
		&= \ve_i\tran\Big[\alpha_j\conj\ve_j +\ \sum_{\mathclap{k\in\aOut{\sBE}(j)}}\aOut{\beta_k\conj}\ve_{\fbE(k)} +\ \sum_{\mathclap{k\in\aIn{\sBE}(j)}}\aIn{\beta_k\conj}\ve_{\faE(k)}\Big]\allowdisplaybreaks\\
		&= \left\{
			\begin{array}{ll}
				\alpha_j\conj & \text{if}\ i=j\\
				\aIn{\beta_{k_1}}\conj & \text{if}\ \exists\,k_1\in\sE:\faE(k_1)=i\;\land\;\fbE(k_1)=j\\
				\aOut{\beta_{k_2}}\conj & \text{if}\ \exists\,k_2\in\sE:\faE(k_2)=j\;\land\;\fbE(k_2)=i\\
				0 & \text{otherwise}\;.
			\end{array}
		\right.
\end{align*}
Analogously, this analysis can be performed for the matrices $\mIout_k$, $\mIin_k$, $\aLB{\mM_k}$, $\aUB{\mM_k}$, $\mA_k$, $\aLB{\mA_k}$, and $\mAub_k$, where the corresponding result is documented in Table~\ref{tab:acbmtxelements}. It follows that all off-diagonal elements of $\mM_n$, $-\mM_n$, $\aLB{\mM_k}$, and $\aUB{\mM_k}$ are zero, for all $n\in\sV$ and $k\in\sE$. For the other matrices, it can be observed that their off-diagonal element in row $i$ and column $j$, with $i\neq j$, may be nonzero if and only if there exists an AC branch between bus $i$ and bus $j$. Therefore, if there does not exist any AC branch between bus $i$ and $j$, then $[\fPsiCone(\vLambda)]_{i,j}=0$. To complete the proof, assume there exists an AC branch $k_1\in\sE$ with source bus $i\in\sV$ and destination bus $j\in\sV$, i.e., $\exists\,k_1\in\sE:\faE(k_1)=i\;\land\;\fbE(k_1)=j$. Then,
$[\mP_j]_{i,j}=\aIn{\beta_{k_1}}\conj/2$,
$[\mP_i]_{i,j}=\aOut{\beta_{k_1}}/2$,
$[\mQ_j]_{i,j}=\aIn{\beta_{k_1}}\conj/(2\iu)$,
$[\mQ_i]_{i,j}=-\aOut{\beta_{k_1}}/(2\iu)$,
$[\mIout_{k_1}]_{i,j}=\aOut{\alpha_{k_1}\conj}\aOut{\beta_{k_1}}$,
$[\mIin_{k_1}]_{i,j}=\aIn{\alpha_{k_1}}\aIn{\beta_{k_1}\conj}$,
$[\mA_{k_1}]_{i,j}=-1$,
$[\aLB{\mA_{k_1}}]_{i,j}=\tan(\aLB{\delta_k}) + \iu$, and
$[\mAub_{k_1}]_{i,j}=-\tan(\aUB{\delta_k}) - \iu$,
whereas the element in row $i$ and column $j$ is zero for all other matrices. Consequently,
\begin{align*}
	&[\fPsiCone(\vLambda)]_{i,j}\in\sC_{k_1}\quad\text{if}\ 
	\exists\,k_1\in\sE:\faE(k_1)=i\;\land\;\fbE(k_1)=j\;.
\intertext{Analogously, in case of an AC branch in the opposite direction, it can be concluded from the conjugate symmetry of the constraint matrices that}
	&[\fPsiCone(\vLambda)]_{i,j}\in\sC_{k_2}\conj\quad\text{if}\ 
	\exists\,k_2\in\sE:\faE(k_2)=j\;\land\;\fbE(k_2)=i
\end{align*}
which completes the proof.

\section{Proof of Theorem~\ref{thm:crtconeishalfplane}}
	\label{apx:thm:crtconeishalfplane}

Corollary~\ref{cor:convhullarg} states that the conic hull $\sC_k = \cone(\sX_k)$ can be described via the arguments of the elements of $\sX_k$. To this end, let $\varphi_k=\arg(-\rho_k)=\pi+\arg(\rho_k)$ to express the arguments of the elements of $\sX_k$ as follows.
\begin{align*}
	\arg(\aOut{\beta_{k}}/2)
		&=\arg(-\rho_k\aBar{y_k}/2)
		 =\varphi_k+\arg(\aBar{y_k})\allowdisplaybreaks\\
	\arg(\aIn{\beta_{k}}\conj/2)
		&=\arg(-\rho_k\aBar{y_k}\conj/2)
		 =\varphi_k-\arg(\aBar{y_k})\allowdisplaybreaks\\
	\arg(-\aOut{\beta_{k}}/(2\iu))
		&=\arg(\rho_k\aBar{y_k}/(2\iu))
		 =\varphi_k+\pi/2+\arg(\aBar{y_k})\allowdisplaybreaks\\
	\arg(\aIn{\beta_{k}}\conj/(2\iu))
		&=\arg(-\rho_k\aBar{y_k}\conj/(2\iu))
		 =\varphi_k-\pi/2-\arg(\aBar{y_k})\allowdisplaybreaks\\
	\arg(\aOut{\alpha_{k}\conj}\aOut{\beta_{k}})
		&=\arg(\abs{\aOut{\rho_k}}^2(\aBar{y_k}\conj+\aOut{y_k}\conj)(-\rho_k\aBar{y_k}))\\
		&=\arg(-\rho_k\abs{\aOut{\rho_k}}^2\abs{\aBar{y_k}}^2[1+\aOut{y_k}\conj\aBar{y_k}/\abs{\aBar{y_k}}^2])\\
		&=\varphi_k+\arg(1+\aOut{y_k}\conj\aBar{y_k}/\abs{\aBar{y_k}}^2)\allowdisplaybreaks\\
	\arg(\aIn{\alpha_{k}}\aIn{\beta_{k}\conj})
		&=\arg(\abs{\aIn{\rho_k}}^2(\aBar{y_k}+\aIn{y_k})(-\rho_k\aBar{y_k}\conj))\\
		&=\arg(-\rho_k\abs{\aIn{\rho_k}}^2\abs{\aBar{y_k}}^2[1+\aIn{y_k}\aBar{y_k}\conj/\abs{\aBar{y_k}}^2])\\
		&=\varphi_k+\arg(1+\aIn{y_k}\aBar{y_k}\conj/\abs{\aBar{y_k}}^2)\allowdisplaybreaks\\
	\arg(-1)
		&=\pi
		 =\varphi_k-\arg(\rho_k)\allowdisplaybreaks\\
	\arg(\tan(\aLB{\delta_k}) &+ \iu)
		 =\pi/2-\aLB{\delta_k}\allowdisplaybreaks\\
	\arg(-\tan(\aUB{\delta_k}) &- \iu)
		 =3\pi/2-\aUB{\delta_k}
\end{align*}
Definition~\ref{def:passiveacbranches} and~\ref{def:inductiveseriesadmittance} imply that $\aBar{y_k}$ is in the fourth quadrant of the $\nC$-plane, thus $-\pi/2\leq\arg(\aBar{y_k})\leq 0$. From Definition~\ref{def:properinsulation} it follows that
\begin{align}
	\abs{\aOut{y_k}\conj\aBar{y_k}/\abs{\aBar{y_k}}^2}
		&= \abs{\aOut{y_k}}\abs{\aBar{y_k}}/\abs{\aBar{y_k}}^2
		 = \abs{\aOut{y_k}}/\abs{\aBar{y_k}} \leq 1\\
	\abs{\aIn{y_k}\aBar{y_k}\conj/\abs{\aBar{y_k}}^2}
		&= \abs{\aIn{y_k}}\abs{\aBar{y_k}}/\abs{\aBar{y_k}}^2
		 = \abs{\aIn{y_k}}/\abs{\aBar{y_k}} \leq 1\;.
\end{align}
Therefore, the values $1+\aOut{y_k}\conj\aBar{y_k}/\abs{\aBar{y_k}}^2$ and $1+\aIn{y_k}\aBar{y_k}\conj/\abs{\aBar{y_k}}^2$ are in the right half of the $\nC$-plane, which implies
\begin{align}
	-\pi/2 &\leq \arg(1+\aOut{y_k}\conj\aBar{y_k}/\abs{\aBar{y_k}}^2) \leq \pi/2\\
	-\pi/2 &\leq \arg(1+\aIn{y_k}\aBar{y_k}\conj/\abs{\aBar{y_k}}^2) \leq \pi/2\;.
\end{align}
Furthermore, Definition~\ref{def:totalphaseshift} states that $-\pi/2\leq\arg(\rho_k)\leq\pi/2$ and from~\eqref{eqn:sysmodel:deltabndcrt} it follows that
\begin{align}
	\varphi_k-\pi/2 \leq \pi/2&-\aLB{\delta_k} < \varphi_k-\arg(\rho_k)
	\\
	\varphi_k-\arg(\rho_k) < 3&\pi/2-\aUB{\delta_k} \leq \varphi_k+\pi/2\;.
\end{align}
Consequently, for all $x\in\sX_k$, $\arg(x)\in[\varphi_k-\pi/2,\varphi_k+\pi/2]$. In conjunction with Corollary~\ref{cor:convhullarg} this implies that
\begin{equation}
	\sC_k \subseteq \sH_k=\cone\big(\en^{\iu(\varphi_k-\pi/2)},\en^{\iu\varphi_k},\en^{\iu(\varphi_k+\pi/2)}\big)\;.
\end{equation}
It can be observed that $\sH_k\subset\nC$ is a half-space defined by the normal $p=\en^{\iu(\varphi_k-\pi)}=\en^{\iu\arg(\rho_k)}$ and, as a consequence of Proposition~\ref{prp:scalednormal}, $\rho_k=\abs{\rho_k}\en^{\iu\arg(\rho_k)}$ is also a normal of $\sH_k$, which completes the proof.
\begin{remark}
Note that the exclusion of lower bounds on the power injection in~\eqref{eqn:opfoptprobfull} originates from the proof above. They introduce $\pi$-rotated values that, in general, extend the cones $\sC_k$ to the entire complex plane and invalidate Theorem~\ref{thm:crtconeishalfplane}.
\end{remark}

\section{Proof of Theorem~\ref{thm:objmtxelements}}
	\label{apx:thm:objmtxelements}

Due to the parametrization for the OPF problem,
\begin{equation}
	[\mC_0]_{i,j} = [\mC]_{i,j} = [w\aCost{\mC}]_{i,j} + [\aLoss{\gamma}\aLoss{\mC}]_{i,j}\;.
		\label{eqn:c0elements}
\end{equation}
For the first summand, it follows from~\eqref{eqn:costobjparamdef} that
\begin{equation}
	[w\aCost{\mC}]_{i,j} = \sum_{n\in\sV} w\gamma_n[\mP_n]_{i,j}
\end{equation}
and as $w\gamma_n\in\nRnn$, for all $n\in\sV$, this is equivalent to $[\fPsiCone(\vLambda)]_{i,j}$ in~\eqref{eqn:kktpsiconeelements} with a corresponding choice of $\vLambda$. Thus, it follows from Theorem~\ref{thm:crtmtxcone} and~\ref{thm:crtconeishalfplane} that, for $i,j\in\sV$ and $i\neq j$,
\begin{align*}
	[w\aCost{\mC}]_{i,j}
		&\in\left\{
			\begin{array}{ll}
				\sH_{k_1}
					& \text{if}\ \exists\,k_1\in\sE:\faE(k_1)=i\;\land\;\fbE(k_1)=j\\
				\sH_{k_2}\conj
					& \text{if}\ \exists\,k_2\in\sE:\faE(k_2)=j\;\land\;\fbE(k_2)=i\\
				\{0\} & \text{otherwise}\;.
			\end{array}
		\right.
\end{align*}
For the second summand in~\eqref{eqn:c0elements}, it follows from~\eqref{eqn:lossobjparamdef} that
\begin{align}
	[\aLoss{\gamma}\aLoss{\mC}]_{i,j} = \sum_{k\in\sE}\aLoss{\gamma}[\aLoss{\mP_k}]_{i,j} + \sum_{n\in\sV}\aLoss{\gamma}\real(\aShunt{y_n})[\mM_n]_{i,j}\;.
\end{align}
For $i\neq j$, the second term is zero as $\mM_n$ is diagonal for all $n\in\sV$. For the first term, it follows from the definition of $\aLoss{\mP_k}$ that $[\aLoss{\mP_k}]_{i,j} = ([\aLoss{\mS_k}]_{i,j}+[\aLoss{\mS_k}]_{j,i}\conj)/2$, where it can be shown that $[\aLoss{\mS_k}]_{i,j}=\ve_i\tran\aLoss{\mS_k}\ve_j$ comprises
\begin{equation}
	[\aLoss{\mS_k}]_{i,j}
		= \left\{
			\begin{array}{ll}
			\aOut{\alpha_k\conj} & \text{if}\ i=j=\faE(k)\\
			\aIn{\alpha_k\conj}  & \text{if}\ i=j=\fbE(k)\\
			\aIn{\beta_k\conj}   & \text{if}\ i=\faE(k)\;\land\;j=\fbE(k)\\
			\aOut{\beta_k\conj}  & \text{if}\ i=\fbE(k)\;\land\;j=\faE(k)\\
			0                    & \text{otherwise}\;.
			\end{array}
			\right.
\end{equation}
Therewith, $[\aLoss{\mP_k}]_{i,j}$ can be specified as
\begin{align}
	[\aLoss{\mP_k}]_{i,j}
		= \left\{
			\begin{array}{ll}
			\real(\aOut{\alpha_k})			& \text{if}\ i=j=\faE(k)\\
			\real(\aIn{\alpha_k})			& \text{if}\ i=j=\fbE(k)\\
			-\rho_k\real(\aBar{y_k})	& \text{if}\ i=\faE(k)\;\land\;j=\fbE(k)\\
			-\rho_k\conj\real(\aBar{y_k}) 		& \text{if}\ i=\fbE(k)\;\land\;j=\faE(k)\\
			0								& \text{otherwise}
			\end{array}
			\right.
\end{align}
in which it is recognized that $(\aOut{\beta_k}+\aIn{\beta_k\conj})/2 = -\rho_k\real(\aBar{y_k})$. From Corollary~\ref{cor:seriesadmittance} it follows that
\begin{equation}
	\real(\rho_k\conj [-\rho_k\real(\aBar{y_k})])
		= -\abs{\rho_k}^2\real(\aBar{y_k}) < 0
\end{equation}
and, with $\sH_k$ in~\eqref{eqn:hndef}, $[\aLoss{\mP_k}]_{i,j}$ with $i\neq j$ satisfies
\begin{align*}
	[\aLoss{\mP_k}]_{i,j}
		&\in\left\{
			\begin{array}{ll}
				\!\interior(\sH_{k_1})
					& \text{if}\ \exists\,k_1\in\sE:\faE(k_1)=i\;\land\;\fbE(k_1)=j\\
				\!\interior(\sH_{k_2}\conj)
					& \text{if}\ \exists\,k_2\in\sE:\faE(k_2)=j\;\land\;\fbE(k_2)=i\\
				\!\{0\} & \text{otherwise}\;.
			\end{array}
		\right.
\end{align*}
With Definition~\ref{def:noparallelacbranches} and~\ref{def:lossweightnonzero}, this implies that
\begin{align*}
	[\aLoss{\gamma}\aLoss{\mC}]_{i,j}
		&\in\left\{
			\begin{array}{ll}
				\!\interior(\sH_{k_1})
					& \!\!\text{if}\ \exists\,k_1\in\sE:\faE(k_1)=i\;\land\;\fbE(k_1)=j\\
				\!\interior(\sH_{k_2}\conj)
					& \!\!\text{if}\ \exists\,k_2\in\sE:\faE(k_2)=j\;\land\;\fbE(k_2)=i\\
				\!\{0\} & \!\!\text{otherwise}
			\end{array}
		\right.
\end{align*}
for $i,j\in\sV$ and $i\neq j$. With this characterization of the summands in~\eqref{eqn:c0elements}, Lemma~\ref{lem:halfplanesuminint} completes the proof.

\section*{Acknowledgment}

The authors would like to thank Prof. Thomas Hamacher of Technische Universit\"at M\"unchen (TUM) for his support and M.~Hotz would like to thank Matthias Huber and Dominic Hewes of TUM for the valuable discussions.

\IEEEtriggeratref{39}

\bibliographystyle{IEEEtran}
\bibliography{hotzbib}

\end{document}